\input ./style/arxiv-general.cfg
\documentclass[MSNbibl,number,citesort,secthm,seceqn,dvips]{arxbj}
\makeatletter
   \@ifpackageloaded{graphicx}{}{\usepackage{graphicx}}
\makeatother
\usepackage{upgreek}
\usepackage{mathrsfs}

\volume{22}
\issue{3}
\pubyear{2016}
\firstpage{1671}
\lastpage{1708}
\doi{10.3150/15-BEJ707}% Updated by VTEXPTS2LaTeX.exe, 04.05.2015 15:10
\docsubty{FLA}

\makeatletter

\newcommand{\eqref}[1]{(\ref{#1})}

\newtheorem{lem}[thm]{Lemma}
\newtheorem{prop}[thm]{Proposition}
\newremark{rem}[thm]{Remark}
\newremark{exm}[thm]{Example}
\newproclaim{asm}[thm]{Assumption}
\newproclaim{defn}[thm]{Definition}
\makeatother

\begin{document}
\begin{frontmatter}

\title{Functional limit theorems for generalized variations of the
fractional Brownian sheet}
\runtitle{Limit theorems for the fractional Brownian sheet}

\begin{aug}
%%%% inicialai - be tarpu
% Corresponding author: Mikko Pakkanen - m.pakkanen@imperial.ac.uk% Updated by VTEXPTS2LaTeX.exe, 05.05.2015 07:23
%Updated by VTEXPTS2LaTeX.exe, 04.05.2015 15:10
\author[A]{\inits{M.S.}\fnms{Mikko S.}~\snm{Pakkanen}\corref{}\thanksref{A,B}\ead[label=e1]{m.pakkanen@imperial.ac.uk}}
\and
\author[C]{\inits{A.}\fnms{Anthony}~\snm{R\'eveillac}\thanksref{C}\ead[label=e2]{anthony.reveillac@insa-toulouse.fr}}
%%\runauthor{} %% auto
%\dedicated{}
\address[A]{Department of Mathematics, Imperial College London, South
Kensington Campus, London SW7 2AZ, UK. \printead{e1}}
\address[B]{CREATES, Aarhus University, Denmark}
\address[C]{INSA de Toulouse, IMT UMR CNRS 5219, Universit\'e de
Toulouse, 135 avenue de Rangueil, 31077 Toulouse Cedex 4, France. \printead{e2}}
\end{aug}

% HISTORY:
%
\received{\smonth{8} \syear{2014}}% Updated by VTEXPTS2LaTeX.exe,
%04.05.2015 15:10
%
\revised{\smonth{1} \syear{2015}}% Updated by VTEXPTS2LaTeX.exe,
%04.05.2015 15:10

% ABSTRACT
\begin{abstract}
We prove functional central and non-central limit theorems for
generalized variations of the anisotropic $d$-parameter fractional
Brownian sheet (fBs) for any natural number $d$. Whether the central or
the non-central limit theorem applies depends on the Hermite rank of
the variation functional and on the smallest component of the Hurst
parameter vector of the fBs. The limiting process in the former result
is another fBs, independent of the original fBs, whereas the limit
given by the latter result is an Hermite sheet, which is driven by the
same white noise as the original fBs. As an application, we derive
functional limit theorems for power variations of the fBs and discuss
what is a proper way to interpolate them to ensure functional convergence.
\end{abstract}

% KEYWORDS
% visi is mazosios raides ir pagal abecele
\end{frontmatter}

\section{Introduction}

Since the seminal works by Breuer and Major \cite{BM}, Dobrushin and
Major \cite{DM}, Giraitis and Surgailis \cite{GS}, Rosenblatt \cite
{Rosenblatt} and Taqqu \cite{Taqqu75,Taqqu1977,Taqqu78,Taqqu79}, much
attention has been given to the study of the asymptotic behaviour of
normalized functionals of Gaussian fields, as these quantities arise
naturally in applications, for example, where models exhibiting \emph
{long-range dependence} are needed. The aforementioned papers focus on
nonlinear functionals of a stationary Gaussian field, for which one can
derive a \emph{central limit theorem} (in a finite-dimensional sense or
in a functional sense) if the correlation function of the field decays
sufficiently fast to zero; see \cite{BM} for a precise formulation.
However, if the correlation function decays too slowly to zero, then
only a \emph{non-central limit theorem} can be established, meaning
that the limiting distribution fails to be Gaussian; see, for example,
\cite{Rosenblatt}.

In particular, these results apply to functionals of the \textit
{fractional Brownian motion} (fBm). Let $B_H:=\{B_H(t) \dvt  t\in\mathbb{R}\}
$ be a fBm with Hurst parameter $H\in(0,1)$, which is the unique (in
law) $H$-self similar Gaussian process with stationary increments; see
\eqref{ssi} and \eqref{si} below for the definitions of these key properties.
The behaviour of the so-called \textit{Hermite variations} of $B_H$,
depending on the value of $H$, can be described as follows. Let $k \in
\{1,2,\ldots\}$ and let $P_k$ denote the $k$th Hermite polynomial, the
definition of which we recall in \eqref{Hermite-poly} below. Applying
results from \cite{BM,DM,GS,Taqqu79}, one obtains that
\begin{longlist}[(a)]
\item[(a)] If $H\in (0,1-\frac{1}{2k} )$, then
\begin{eqnarray*}
&& n^{-1/2} \sum_{j=1}^n
P_k \biggl(n^H \biggl(B_H \biggl(
\frac
{j}{n} \biggr)-B_H \biggl(\frac{j-1}{n} \biggr) \biggr)
\biggr) \mathop{\longrightarrow}_{n\rightarrow\infty}^{\mathscr{L}} N \bigl(0,
\sigma^2_1(H,k) \bigr).
\end{eqnarray*}
\item[(b)] If $H=1-\frac{1}{2k}$, then
\begin{eqnarray*}
&& \bigl(n \log(n) \bigr)^{-1/2} \sum_{j=1}^n
P_k \biggl(n^H \biggl(  {B_H \biggl(
\frac{j}{n} \biggr)-B_H \biggl(\frac{j-1}{n} \biggr)}
\biggr) \biggr) \mathop{\longrightarrow}_{n\rightarrow\infty}^{\mathscr{L}} N\biggl(0,
\sigma ^2_1\biggl({ 1-\frac{1}{2k}},k
\biggr)\biggr).
\end{eqnarray*}
\item[(c)] If $H \in (1-\frac{1}{2k},1 )$, then
\begin{eqnarray*}
&& n^{1-2H} \sum_{j=1}^n
P_k \biggl(n^H \biggl( {B_H \biggl(
\frac
{j}{n} \biggr)-B_H \biggl(\frac{j-1}{n} \biggr)}
\biggr) \biggr) \mathop{\longrightarrow}_{n\rightarrow \infty
}^{L^2(\Omega)}
\operatorname{Hermite}_{1,k} \bigl(1-k(1-H) \bigr).
\end{eqnarray*}
\end{longlist}
Above, $\displaystyle{\mathop{\rightarrow}^{\mathscr{L}}}$ denotes convergence in law,
$N
(0,\sigma^2_1(H,k) )$ denotes the centered Gaussian law with
variance $\sigma^2_1(H,k)>0$, whereas $\operatorname{Hermite}_{1,k}
(1-k(1-H) )$ stands for a so-called Hermite random variable given by
the value of an \emph{Hermite process}, of order $k$ with Hurst
parameter $1-k(1-H)\in(\frac{1}{2},1)$, at time $1$. Such an Hermite
process can be represented as a $k$-fold multiple Wiener integral with
respect to Brownian motion, as proven by Taqqu \cite{Taqqu78,Taqqu79}.
Moreover, the process is non-Gaussian if $k\geq2$. (More details on
the Hermite process are provided in Section~\ref{ssecmain}.)
The key observation here is that there are two regimes: \emph
{Gaussian}, subsuming cases (a) and (b), and \emph
{Hermite}, case (c), depending on the Hurst parameter $H$ and on
the order $k$.

The convergences in all cases (a), (b), and (c)
can be extended to more general functionals, which we call \emph
{generalized variations} in this paper, obtained by replacing the
Hermite polynomial $P_k$ with a function
\begin{equation}
\label{introf}
f(u):=\sum_{k=\underline{k}}^\infty
a_k P_k(u),\qquad   u \in\mathbb{R},
\end{equation}
where $\underline{k}$ is the so-called \textit{Hermite rank} of $f$.
(Naturally, conditions on the summability of the coefficients
$a_{\underline{k}},a_{\underline{k}+1},\ldots$ have to be added.) In
this setting, the prevailing regime (Gaussian or Hermite) will depend
on the Hurst parameter $H$ and on the Hermite rank $\underline{k}$
analogously to the simpler setting discussed above.
In addition, functional versions of these asymptotic results (under
additional assumptions on the coefficients $a_{\underline
{k}},a_{\underline{k}+1},\ldots$) can be proven in the Skorohod space
$D([0,1])$; see \cite{Taqqu75,Taqqu79}.

In connection to applications that involve spatial or spatio-temporal
modeling, processes of multiple parameters are also of interest.
Recently, there has been interest in understanding the asymptotic
behaviour of realized quadratic variations and power variations of
\emph
{ambit fields} \cite{BG2011,P13}. An ambit field is an anisotropic
multiparameter process driven by white noise, or more generally, by an
infinitely-divisible random measure. The problem of finding
distributional limits (central or non-central limit theorems) for such
power variations is, however, intricate because the dependence
structure of an ambit field can be very general; only a ``partial''
central limit theorem is obtained in \cite{P13}. As a first
approximation, it is thus useful to study this problem with simpler
processes that incorporate some of the salient features of ambit
fields, such as the non-semimartingality of one-parameter ``marginal
processes'' (see \cite{P13}, Section~2.2) and strong dependence.
A tractable process that incorporates some key features of ambit fields
is the \emph{fractional Brownian sheet} (fBs), defined by Ayache \textit{et al.} \cite{ALP2002}, which is a multi-parameter extension of the fBm.
In particular, it is a Gaussian process with stationary rectangular increments.

For concreteness, let $Z:=\{ Z(t) \dvt t \in[0,1]^2 \}$ be a
two-parameter anisotropic fBs with Hurst parameter $(H_1,H_2) \in
(0,1)^2$; see Section~\ref{ssecfbs} for a precise definition. In view
of the asymptotic behaviour in cases (a), (b), and~(c) involving the fBm, it is natural to ask what is the asymptotic
behaviour of Hermite variations of $Z$ with different values of $H_1$
and $H_2$. Consider, for example, the ``mixed'' case where $H_1<1-\frac
{1}{2k}$ and $H_2>1-\frac{1}{2k}$, which has no counterpart in the
one-parameter setting. Because of the structure of the fBs, it is
tempting to conjecture that in this case the limiting law is a mixture
of a Gaussian law and a marginal law of an Hermite process. However, as
shown in \cite{RST2012}, this is not the case and once again only two
limiting laws can be obtained:
\begin{longlist}[(a$'$)]
\item[(a$'$)] If $(H_1,H_2)\in(0,1)^2 \setminus (1-\frac{1}{2k},1
)^2$, then
\begin{eqnarray*}
&& \varphi(n,H_1,H_2) \sum_{j_1=1}^n
\sum_{j_2=1}^n P_k
\biggl(n^{H_1+H_2}Z \biggl( \biggl[\frac{j_1-1}{n},\frac{j_1}{n} \biggr)
\times \biggl[\frac{j_2-1}{n},\frac{j_2}{n} \biggr) \biggr) \biggr)
\\
&&\quad\mathop{\longrightarrow}_{n\rightarrow\infty}^{\mathscr{L}} N \bigl(0,\sigma
^2_2(H_1,H_2,\underline {k})
\bigr).
\end{eqnarray*}
\item[(b$'$)] If $(H_1,H_2)\in (1-\frac{1}{2k},1 )^2$, then
\begin{eqnarray*}
&& \varphi(n,H_1,H_2) \sum_{j_1=1}^n
\sum_{j_2=1}^n P_k
\biggl(n^{H_1+H_2}Z \biggl( \biggl[\frac{j_1-1}{n},\frac{j_1}{n} \biggr)
\times \biggl[\frac{j_2-1}{n},\frac{j_2}{n} \biggr) \biggr) \biggr)
\\
&&\quad\mathop{\longrightarrow}_{n\rightarrow\infty}^{L^2(\Omega)} \operatorname{Hermite}_{2,k}
\bigl(1-k(1-H_1),1-k(1-H_2) \bigr).
\end{eqnarray*}
\end{longlist}
Above, $Z ( [\frac{j_1-1}{n},\frac{j_1}{n} ) \times
[\frac
{j_2-1}{n},\frac{j_2}{n} )  )$ stands for the increment of $Z$
over the rectangle $ [\frac{j_1-1}{n},\frac{j_1}{n} ) \times
[\frac{j_2-1}{n},\frac{j_2}{n} )$, defined in Section~\ref{ssecincrvar} below, and $\varphi(n,H_1,H_2)$ is a suitable scaling
factor; see \cite{RST2012}, pages 9--10, for its definition. The limit in
the case (b$'$) is the value of a two-parameter \emph{Hermite
sheet} (see Section~\ref{ssecmain}), of order $k$ with Hurst parameter
$ (1-k(1-H_1),1-k(1-H_2) )\in(\frac{1}{2},1)^2$, at point
$(1,1)$. Contrary to the one-parameter case, the results obtained in
\cite{RST2012} are proved only for one-dimensional laws; neither
finite-dimensional (except in the particular setting of \cite
{Rev2009b}) nor functional convergence (i.e., tightness in a function
space) of Hermite variations has been established so far. (In
particular in the $d$-parameter realm with $d\geq2$, tightness is a
non-trivial issue, which has not been addressed in \cite{RST2012} or in
the related paper \cite{Rev2009b}.)

The first main result of this paper addresses the question about
functional convergence in the general, $d$-parameter case for any $d
\in\mathbb{N}$.
We prove a functional central limit theorem, Theorem~\ref{gvar-clt},
for generalized variations of a $d$-parameter anisotropic fBs $Z$.
(As mentioned above, generalized variations extend Hermite variations
by replacing $P_k$ with a function $f$ of the form \eqref{introf}.)
This result applies if at least one of the components of the Hurst
parameter vector $H=(H_1,\ldots,H_d)\in(0,1)^d$ of $Z$ is less than or
equal to $1-\frac{1}{2\underline{k}}$, where $\underline{k}$ is the
Hermite rank of $f$. A novel feature of this result is that the
limiting process is a new fBs, independent of $Z$, with Hurst parameter
vector $\widetilde{H} =  (\widetilde{H}_1,\ldots,\widetilde
{H}_d
)$ given by
\begin{eqnarray*}
&&\widetilde{H}_\nu:= %
\cases{ \displaystyle\frac{1}{2}, &\quad
$\displaystyle H_\nu\leq1-\frac{1}{2\underline{k}}$,\vspace*{3pt}
\cr
1-\underline{k}(1-H_\nu), &\quad $\displaystyle H_\nu> 1-\frac{1}{2\underline{k}}$,
}
\end{eqnarray*}
for $\nu\in\{1,\ldots,d\}$.
Note, in particular, that if $H \in (0,1-\frac{1}{2\underline
{k}} ]^d$, then the limit reduces to an ordinary Brownian sheet. The
proof of Theorem~\ref{gvar-clt} is based on the limit theory for
multiple Wiener integrals, due to Nualart and Peccati \cite{NP2005},
and its multivariate extension by Peccati and Tudor \cite{PT2005}. To
prove the functional convergence asserted in Theorem~\ref{gvar-clt}, we
use the tightness criterion of Bickel and Wichura \cite{BW1971} in the
space $D([0,1]^d)$, which is $d$-parameter generalization of
$D([0,1])$, and a moment bound for nonlinear functionals of a
stationary Gaussian process on $\mathbb{Z}^d$ (Lemma~\ref{moment-bound}).

The second main result of this paper is a functional non-central limit
theorem, Theorem~\ref{gvar-nclt}, for generalized variations of $Z$ in
the remaining case where each of the components of $H$ is greater than
$1-\frac{1}{2\underline{k}}$. In this case, the limit is a
$d$-parameter Hermite sheet and the convergence holds in probability
and also pointwise in $L^2(\Omega)$. Assuming that $Z$ is defined by a
moving-average representation with respect to a white noise $\mathscr
{W}$ on $\mathbb{R}^d$, we can give a novel and explicit description
of the
limit; it is defined using the representation introduced by Clarke De
la Cerda and Tudor \cite{CT2012} with respect to the same white noise
$\mathscr{W}$. This makes the relation between $Z$ and the Hermite
sheet precise and constitutes a step further compared to the existing
literature (see \cite{NNT,RST2012}), where the limiting Hermite
process/sheet is simply obtained as an abstract limit of a Cauchy
sequence, from which the properties of the limiting object are deduced.

As an application of Theorems \ref{gvar-clt} and \ref{gvar-nclt}, we
study the asymptotic behaviour of power variations of the fBs $Z$. As a
straightforward consequence of our main results, we obtain a law of
large numbers for these power variations. We then study the more
delicate question regarding the asymptotic behaviour of rescaled
fluctuations of power variations around the limit given by the law of
large numbers. In the case of \emph{odd} power variations, the rescaled
fluctuations have a limit, either Gaussian or Hermite, but with \emph
{even} power variations, the fluctuations might not converge in a
functional way if $d \geq2$. We show that this convergence issue does
not arise at all if one considers instead continuous, multilinear
interpolations of power variations.

The paper is organized as follows. In Section~\ref{secmainresults}, we
introduce the setting of the paper, some key definitions and the
statements of Theorems \ref{gvar-clt} and \ref{gvar-nclt}. The proofs
of these two main results are presented in Sections~\ref{fdd-convergence} and~\ref{fun-convergence},
the former section collecting the finite-dimensional and the latter
the functional arguments. Finally, the application to power variations
is given in Section~\ref{secpower-var}.

\section{Preliminaries and main results}\label{secmainresults}

\subsection{Notation}\label{ssecnotations}

We use the convention that $\mathbb{N}:=\{1,2,\ldots\}$ and $\mathbb{R}_+
:=[0,\infty)$. The notation $|A|$ stands for the cardinality of a
finite set $A$. For any $y \in\mathbb{R}$, we write $\lfloor y
\rfloor
:=
\max\{n \in\mathbb{Z}\dvt n \leq v \}$, $\{y\} :=y - \lfloor y \rfloor$,
and $(y)_+ :=\max(y,0)$. The symbol $\gamma$ denotes the standard
Gaussian measure on $\mathbb{R}$, that is, $\gamma(\mathrm{d}y) :=(2\uppi)^{-1/2}
\exp(-y^2/2) \,\mathrm{d}y$. From now on we fix $d$ in $\mathbb{N}$.

For any vectors $s = (s_1,\ldots,s_d) \in\mathbb{R}^d$ and $t =
(t_1,\ldots
,t_d) \in\mathbb{R}^d$, the relation $s \leq t$ (resp., $s < t$) signifies
that $s_\nu\leq t_\nu$ (resp., $s_\nu< t_\nu$) for all $\nu\in\{
1,\ldots,d\}$.
We also use the notation
\begin{eqnarray*}
st & :=& (s_1t_1,\ldots,s_dt_d)
\in\mathbb{R}^d, \qquad \frac{s}{t} := \biggl(\frac{s_1}{t_1},
\ldots,\frac{s_d}{t_d} \biggr)\in\mathbb{R}^d,
\\
\lfloor s \rfloor& :=& \bigl(\lfloor s_1 \rfloor,\ldots,\lfloor
s_d \rfloor\bigr)\in\mathbb{Z}^d, \qquad \langle s\rangle
:=s_1\cdots s_d \in \mathbb{R} ,
\\
|s| & :=& \bigl(|s_1|,\ldots,|s_d|\bigr) \in\mathbb{R}^d_+,
\qquad \{ s \}  :=\bigl(\{ s_1\},\ldots,\{s_d\}\bigr)
\in[0,1)^d.
\end{eqnarray*}
Further, when $s \in\mathbb{R}^d_+$, we write $s^t
:=(s_1^{t_1},\ldots
,s_d^{t_d}) \in\mathbb{R}^d_+$, and when $s\leq t$, we write $[s,t) :=
[s_1,t_1)\times\cdots\times[s_d,t_d)\subset\mathbb{R}^d$.
Occasionally, we
use the norm $\| s\|_\infty:=\max(|s_1|,\ldots,|s_d|)$ for $s
\in\mathbb{R}^d$.

For the sake of clarity, we will consistently use the following
convention: $i,i^{(1)},i^{(2)},\ldots$ are multi-indices (vectors) in
$\mathbb{Z}^d$ and $j,j_1,j_2,\ldots$ are indices (scalars) in
$\mathbb{Z}$.

\subsection{Anisotropic fractional Brownian sheet}\label{ssecfbs}

We consider an anisotropic, $d$-parameter fractional Brownian sheet
(fBs) $Z:=\{Z(t) \dvt t \in\mathbb{R}^d\}$ with Hurst parameter $H \in
(0,1)^d$, which is a centered Gaussian process with covariance
\begin{equation}
\label{fBscov}
R^{(d)}_{H}(s,t):=\mathbf{E}\bigl[Z(s)Z(t)
\bigr] = \prod_{\nu=1}^d
R^{(1)}_{H_\nu
}(s_\nu,t_\nu), \qquad  s,
  t \in\mathbb{R}^d,
\end{equation}
where
\begin{eqnarray*}
&& R^{(1)}_{H_\nu}(s_\nu,t_\nu) :=
\tfrac{1}{2} \bigl(|s_\nu |^{2H_\nu
}+|t_\nu|^{2H_\nu}-|s_\nu-t_\nu|^{2H_\nu}
\bigr), \qquad  s_\nu,   t_\nu \in\mathbb{R},
\end{eqnarray*}
is the covariance of a fractional Brownian motion with Hurst parameter
$H_\nu$.

In what follows, it will be convenient to assume that the fBs $Z$ has a
particular representation. To this end, let us denote by $\mathscr
{B}_0(\mathbb{R}^d)$ the family of Borel sets of $\mathbb{R}^d$ with
finite Lebesgue measure.
Let $(\Omega,\mathscr{F},\mathbf{P})$ be a complete probability
space that
supports a \emph{white noise} $\mathscr{W} :=\{\mathscr{W}(A)
\dvt A
\in\mathscr{B}_0(\mathbb{R}^d)\}$, which is a centered Gaussian
process with
covariance
\begin{eqnarray*}
&& \mathbf{E}\bigl[\mathscr{W}(A)\mathscr{W}(B)\bigr] = \operatorname{Leb}_d(A
\cap B),  \qquad A,  B \in\mathscr{B}_0\bigl(\mathbb{R}^d
\bigr),
\end{eqnarray*}
where $\operatorname{Leb}_d(\cdot)$ denotes the Lebesgue measure on
$\mathbb{R}^d$.
The process $Z$ can be defined as a \emph{Wiener integral} with respect
to $\mathscr{W}$ (see, e.g., \cite{Nua2006} for the definition), namely
\begin{equation}
\label{fBs-Wiener} Z(t) :=\int G^{(d)}_{H}(t,u) \mathscr{W}(
\mathrm{d}u),  \qquad t \in \mathbb{R}^d,
\end{equation}
where the kernel
\begin{equation}
\label{tensor-kernel}
G^{(d)}_{H}(t,u) :=\prod
_{\nu=1}^d G^{(1)}_{H_\nu}(t_\nu
,u_\nu), \qquad  t,   u \in\mathbb{R}^d,
\end{equation}
is defined using the one-dimensional Mandelbrot--Van Ness \cite
{MVN1968} kernel
\begin{equation}
\label{MVN}
G^{(1)}_{H_\nu}(t_\nu,u_\nu)
:=\frac{1}{\chi(H_\nu)} \bigl((t_\nu- u_\nu)_+^{H_\nu- {1}/{2}}
- (-u_\nu)_+^{H_\nu- {1}/{2}} \bigr),  \qquad  t_\nu,
u_\nu\in\mathbb{R},
\end{equation}
with the normalizing constant
\begin{eqnarray*}
&& \chi(H_\nu) := \biggl(\frac{1}{2H_\nu}+\int_0^\infty
\bigl((1+y)^{H_\nu-{1}/{2}} - y^{H_\nu-{1}/{2}} \bigr) \,\mathrm{d}y
\biggr)^{{1}/{2}}.
\end{eqnarray*}

We refer to \cite{ALP2002} for a proof that the process $Z$ defined via
\eqref{fBs-Wiener} does indeed have the covariance structure \eqref{fBscov}.
The fBs admits a continuous modification (see \cite{BF2005}, page~1040),
so we may assume from now on that $Z$ is continuous.

\subsection{Increments and generalized variations}\label{ssecincrvar}

Given\vspace*{1pt} a function (or a realization of a stochastic process) $h \dvtx
\mathbb{R}^d
\rightarrow\mathbb{R}$, we define the increment of $h$ over the half-open
hyperrectangle $[s,t)\subset\mathbb{R}^d$ for any $s \leq t$ by
\begin{eqnarray*}
&& h\bigl([s,t)\bigr) :=\sum_{i \in\{0,1\}^d}
(-1)^{d-\sum_{\nu=1}^d i_\nu
}h \bigl((1-i)s+it \bigr).
\end{eqnarray*}
(Note that $i_\nu$ above stands for the $\nu$th component of the
multi-index $i$.)
This definition can be recovered by differencing iteratively with
respect to each of the arguments of the function $h$. Thus, the
increment can be seen as a discrete analogue of the partial derivative
$\frac{\partial^d}{\partial t_1 \cdots\partial t_d}$.

\begin{rem}\label{tensorcase}
It is useful to note that if there exists functions $h_\nu \dvtx \mathbb{R}
\rightarrow\mathbb{R}$, $\nu\in\{1,\ldots,d\}$, such that $h(t) =
h_1(t_1)\cdots h_d(t_d)$ for any $t \in\mathbb{R}^d$, then
\begin{eqnarray*}
&& h\bigl([s,t)\bigr) = \prod_{\nu=1}^d
\bigl(h_\nu(t_\nu)-h_\nu(s_\nu)
\bigr),
\end{eqnarray*}
which is easily verified by induction with respect to $d$ using
iterative differencing.
\end{rem}

Let us fix a sequence $ (m(n) )_{n \in\mathbb{N}} \subset
\mathbb{N}^d$ of
multi-indices with the property
\begin{eqnarray*}
&& \underline{m}(n) :=\min \bigl(m_1(n),\ldots,m_d(n)
\bigr)\mathop{\longrightarrow}_{n \rightarrow\infty} \infty
\end{eqnarray*}
and a function $f \in L^2(\mathbb{R},\gamma)$ such that $\int_\mathbb{R}f(u) \gamma
(\mathrm{d}u) = 0$.
Our aim is to study the asymptotic behaviour of a family $ \{
U^{(n)}_f \dvt n \in\mathbb{N} \}$ of $d$-parameter processes, \emph
{generalized variations} of $Z$, defined by
\begin{eqnarray*}
&& U^{(n)}_f(t) :=\sum_{1 \leq i \leq\lfloor m(n)t \rfloor} f
\biggl( \bigl\langle m(n)^H \bigr\rangle Z \biggl( \biggl[
\frac
{i-1}{m(n)},\frac
{i}{m(n)} \biggr) \biggr) \biggr),\qquad   t
\in[0,1]^d,   n \in \mathbb{N}.
\end{eqnarray*}
In this definition, $ \langle m(n)^H  \rangle= m_1(n)^{H_1}
\cdots m_d(n)^{H_d}$ according to the notation and conventions set
forth in Section~\ref{ssecnotations}.
The realizations of $U^{(n)}_f$ belong to the space $D([0,1]^d)$, which
for $d \geq2$ is a generalization of the space $D([0,1])$ of c\`adl\`ag functions on $[0,1]$. We refer to \cite{BW1971}, page~1662, for the
definition of the space $D([0,1]^d)$. In particular, $C([0,1]^d)
\subset D([0,1]^d)$. We endow $D([0,1]^d)$ with the Skorohod topology
described in \cite{BW1971}, page~1662. Convergence to a continuous
function in this topology is, however, equivalent to uniform
convergence (see, e.g., \cite{P13}, Lemma B.2, for a proof in the case $d=2$).

\subsection{Functional limit theorems for generalized
variations}\label{ssecmain}

We will now formulate two functional limit theorems for the family
$\{U^{(n)}_f \dvt n \in\mathbb{N} \}$ of generalized variations,
defined above.
The class of admissible functions $f$ needs to be restricted somewhat,
however, and the choice of $f$ and the Hurst parameter $H$ of $Z$ will
determine which of the limit theorems applies. Also, we need to rescale
$U^{(n)}_f$ in suitable way that, likewise, depends on both $f$ and $H$.

To this end, recall that the \emph{Hermite polynomials},
\begin{equation}
\label{Hermite-poly}
P_{0} (u) :=1, \qquad  P_k(u)
:=(-1)^k \mathrm{e}^{{u^2}/{2}} \,\frac{\mathrm{d}^k}{\mathrm{d}u^k}
\mathrm{e}^{-{u^2}/{2}}, \qquad  u \in\mathbb{R},   k \in\mathbb{N},
\end{equation}
form a complete orthogonal system in $L^2(\mathbb{R},\gamma)$.
Thus, we may expand $f$ in $L^2(\mathbb{R},\gamma)$ as
\begin{equation}
\label{f-Hermite}
f(u) = \sum_{k=\underline{k}}^\infty
a_k P_k(u),
\end{equation}
where the \emph{Hermite coefficients} $a_{\underline
{k}},a_{\underline
{k}+1},\ldots\in\mathbb{R}$ are such that $a_{\underline{k}}\neq0$ and
\begin{equation}
\label{standard-sum}
\sum_{k = \underline{k}}^\infty k!
a^2_k < \infty.
\end{equation}
The index $\underline{k}$ is called the \emph{Hermite rank} of $f$, and
the proviso $\int_\mathbb{R}f(u) \gamma(\mathrm{d}u) = 0$ ensures
that \mbox{$\underline {k} \geq1$}. We will assume that the Hermite coefficients decay
somewhat faster than what \eqref{standard-sum} entails.

\begin{asm}\label{Hermite-assumption}
The Hermite coefficients $a_{\underline{k}},a_{\underline
{k}+1},\ldots$
of the function $f$ satisfy
\begin{eqnarray*}
&& \sum_{k = \underline{k}}^\infty3^{{k}/{2}} \sqrt{k!}
|a_k| < \infty.
\end{eqnarray*}
\end{asm}

Let us define a sequence $ (c(n) )_{n \in\mathbb{N}}\subset
\mathbb{R}^d_+$ of
rescaling factors by setting for any $\nu\in\{1,\ldots,d\}$ and $n
\in\mathbb{N}$,
\begin{eqnarray*}
&& c_\nu(n) := %
\cases{ \displaystyle m_\nu(n)^{2-2\underline{k}(1-H_\nu)},
& \quad $\displaystyle H_\nu\in \biggl(1-\frac
{1}{2\underline{k}},1 \biggr)$,\vspace*{3pt}
\cr
\displaystyle m_\nu(n) \log \bigl( m_\nu(n) \bigr), &
\quad $\displaystyle H_\nu= 1-\frac
{1}{2\underline{k}}$,\vspace*{3pt}\cr
\displaystyle m_\nu(n), &\quad $\displaystyle H_\nu\in \biggl(0,1-\frac{1}{2\underline{k}}
\biggr)$.}
\end{eqnarray*}

\begin{rem}\label{psidelta-comparison}
Note that $\limsup_{n \rightarrow\infty}\frac{m_\nu(n)}{c_\nu(n)}
<\infty$ and that, in fact, $\lim_{n \rightarrow\infty}\frac{m_\nu
(n)}{c_\nu(n)} = 0$ if $H_\nu\in [1-\frac{1}{2\underline
{k}},1 )$.
\end{rem}

Now we can define a family $ \{\overline{U}_f^{(n)} \dvt n \in
\mathbb{N}\}$
of rescaled generalized variations as
\begin{eqnarray*}
&& \overline{U}_f^{(n)}(t) :=\frac{U_f^{(n)}(t)}{\langle
c(n)\rangle
^{{1}/{2}}}, \qquad  t
\in[0,1]^d,   n \in\mathbb{N}.
\end{eqnarray*}

Our first result is the following functional central limit theorem
(FCLT) for generalized variations. Its proof is carried out in Section~\ref{sseccltproof} and Section~\ref{ssectightness}.

\begin{thm}[(FCLT)]\label{gvar-clt}
Let $f$ be as above such that Assumption~\ref{Hermite-assumption} holds
and suppose that $H \in(0,1)^d \setminus (1-\frac{1}{2\underline
{k}},1 )^d$. Then
\begin{equation}
\label{clt-conv}
\bigl(Z,  \overline{U}_f^{(n)} \bigr)
\mathop{\longrightarrow}_{n \rightarrow \infty }^{\mathscr{L}} \bigl(Z,  \Lambda_{H,f}^{{1}/{2}}
\widetilde {Z} \bigr)\qquad \mbox{in }D\bigl([0,1]^d
\bigr)^2,
\end{equation}
where $\widetilde{Z}$ is a $d$-parameter fBs with Hurst parameter
$\widetilde{H}\in [\frac{1}{2},1 )^d$, independent of $Z$
(defined, possibly, on an extension of $(\Omega,\mathscr{F},\mathbf{P}
)$), and
\begin{equation}
\label{lambdadef}
\Lambda_{H,f} :=\sum_{k=\max(\underline{k},2)}^{\infty}
k! a_k^2 \bigl\langle b^{(k)}\bigr\rangle\in
\mathbb{R}.
\end{equation}
The vectors $\widetilde{H}\in [\frac{1}{2},1 )^d$ and $b^{(k)}
\in\mathbb{R}_+^d$, $k\geq\max(\underline{k},2)$, that appear
above are
defined by setting for any $\nu\in\{1,\ldots,d\}$,
\begin{equation}
\label{H-shift}
\widetilde{H}_\nu:= %
\cases{
\displaystyle\frac{1}{2}, & \quad$\displaystyle H_\nu\in \biggl(0,1-\frac{1}{2\underline{k}} \biggr]$,\vspace*{3pt}\cr
\displaystyle 1-\underline{k}(1-H_\nu), & \quad$\displaystyle H_\nu\in \biggl(1-
\frac{1}{2\underline
{k}},1 \biggr)$,}
\end{equation}
and
\begin{equation}
\label{bdef}
\hspace*{-12pt} b^{(k)}_\nu:=
\cases{\displaystyle
 \sum_{j \in\mathbb{Z}} \biggl( \frac{|j+1|^{2H_\nu}
-2|j|^{2H_\nu}+|j-1|^{2H_\nu}}{2}
\biggr)^{k}, & \hspace*{7pt} $\displaystyle H_\nu\in \biggl(0,1-\frac{1}{2\underline{k}}
\biggr)$,
\vspace*{3pt}\cr
\displaystyle 2 \biggl(\frac{(2\underline{k}-1)(\underline
{k}-1)}{2\underline{k}^2} \biggr)^{\underline{k}} =:\iota (
\underline{k}), &\hspace*{7pt}  $\displaystyle H_\nu= 1-\frac{1}{2\underline{k}},   k = \underline{k}$,\vspace*{3pt}
\cr
\displaystyle\frac{H_\nu^{\underline{k}}(2H_\nu- 1)^{\underline
{k}}}{(1-\underline{k}(1-H_\nu))(1-2\underline{k}(1-H_\nu))} =: \kappa(H_\nu,\underline{k}), &
\hspace*{7pt}$\displaystyle H_\nu\in \biggl(1-\frac
{1}{2\underline
{k}},1 \biggr),   k = \underline{k}$,
\vspace*{3pt}\cr
0, &\hspace*{7pt} $\displaystyle H_\nu\in [1-\frac{1}{2\underline{k}},1 ),   k > \underline{k}$.
}\hspace*{-10pt}
\end{equation}
\end{thm}

\begin{rem}
(1) The counterpart of the convergence \eqref{clt-conv} for finite-dimensional laws holds without Assumption~\ref{Hermite-assumption}; see
Proposition~\ref{fdd-conv-prop} below.

(2) We may use $\max(\underline{k},2)$, instead of $\underline{k}$,
as the lower bound for the summation index $k$ in \eqref{lambdadef}
since $\iota(1) = 0$ and
\begin{eqnarray*}
&& \sum_{j \in\mathbb{Z}} \frac{|j+1|^{2\check{H}} -2|j|^{2\check
{H}}+|j-1|^{2\check{H}}}{2}
\\
&&\quad= \sum_{j \in\mathbb{Z}} \frac{|j|^{2\check
{H}}-|j-1|^{2\check{H}}}{2}-\sum
_{j \in\mathbb{Z}} \frac
{|j|^{2\check
{H}}-|j-1|^{2\check{H}}}{2}= 0
\end{eqnarray*}
for any $\check{H} \in (0,\frac{1}{2} )$. (Then, $\sum_{j
\in\mathbb{Z}
}  |\frac{|j|^{2\check{H}}-|j-1|^{2\check{H}}}{2} | < \infty
$ by
the mean value theorem.)

(3) The convergence \eqref{clt-conv} can be understood in the
framework of \emph{stable convergence} in law, introduced by R\'enyi
\cite{Ren1963}. Equivalently to \eqref{clt-conv}, $\overline
{U}_f^{(n)}$ converges to $\Lambda_{H,f}^{{1}/{2}}\widetilde{Z}$ as
$n \rightarrow\infty$ stably in law with respect to the $\sigma
$-algebra generated by $\{Z(t) \dvt t \in[0,1]^d\}$.
%\item See Appendix~\ref{Gaussian-ma} for an extension of Theorem
%\ref{gvar-clt} beyond the present setting with the fBs.
%
\end{rem}

Theorem~\ref{gvar-clt} excludes the case $H \in (1-\frac
{1}{2\underline{k}},1 )^d$. Then, the generalized variations do have
a limit, but the limit is non-Gaussian, unless $\underline{k}=1$. To
describe the limit, we need the following definition, due to Clarke De
la Cerda and Tudor \cite{CT2012}.

\begin{defn}
An anisotropic, $d$-parameter \emph{Hermite sheet}
$\widehat
{Z} := \{\widehat{Z}(t) \dvt t \in\mathbb{R}^d_+  \}$ of order
$k\geq
2$ with Hurst parameter $\widetilde{H} \in (\frac{1}{2},1 )^d$
is defined as a $k$-fold multiple Wiener integral (see Section~\ref{sseccltproof}) with respect to the white noise $\mathscr{W}$,
\begin{equation}
\label{Hermite-sheet}
\widehat{Z}(t) :=\int\cdots\int\widehat{G}^{(k)}_{\widetilde
{H}}
\bigl(t,u^{(1)},\ldots,u^{(k)} \bigr) \mathscr{W} \bigl(
\mathrm{d} u^{(1)} \bigr)\cdots\mathscr{W} \bigl(\mathrm{d}u^{(k)}
\bigr) :=I^{\mathscr{W}}_k \bigl(\widehat{G}^{(k)}_{\widetilde{H}}
(t,  \cdot  ) \bigr)
\end{equation}
for any $t \in\mathbb{R}^d_+$. In \eqref{Hermite-sheet}, the kernel
$\widehat
{G}^{(d,k)}_{\widetilde{H}}(t,\cdot) \in L^2(\mathbb{R}^{kd})$
is given by
\begin{eqnarray*}
&& \widehat{G}^{(k)}_{\widetilde{H}} \bigl(t,u^{(1)},
\ldots,u^{(k)} \bigr)
\\
&&\quad:=\frac{1}{\widehat{\chi} (\widetilde{H},k )} \int_{[0,t)}\prod
_{\kappa=1}^{k} \prod_{\nu=1}^d
\bigl(y_\nu-u_\nu ^{(\kappa
)} \bigr)_+^{-{1}/{2}-({1-\widetilde{H}_\nu})/{k}}
\,\mathrm{d}y, \qquad  u^{(1)},\ldots,u^{(k)} \in
\mathbb{R}^d,
\end{eqnarray*}
using the normalizing constant
\begin{eqnarray*}
&& \widehat{\chi} (\widetilde{H},k ) :=\prod_{\nu= 1}^d
\biggl(\frac{\mathrm{B} ({1}/{2}-({1-\widetilde{H}_\nu
})/{k},({2(1-\widetilde{H}_\nu)})/{k} )}{\widetilde{H}_\nu(2\widetilde
{H}_\nu-1)} \biggr)^{{1}/{2}},
\end{eqnarray*}
where $\mathrm{B}$ stands for the beta function.
\end{defn}

The Hermite sheet $\widehat{Z}$ is self-similar and has the same
correlation structure as a fBs with Hurst parameter $\widetilde{H}$.
In the case $k=1$, the process $\widehat{Z}$ is Gaussian (in fact, it
coincides with a fractional Brownian sheet with Hurst parameter
$\widetilde{H}$) but for $k \geq2$ it is non-Gaussian. In the case
$k=2$, the name \emph{Rosenblatt sheet} (and \emph{Rosenblatt process},
when $d=1$; see \cite{Tudor08}) is often used, in honor of Murray
Rosenblatt's seminal paper \cite{Ros1961}. See also the recent papers
\cite{MT2013,VT2013} for more details on the Rosenblatt distribution,
including proofs that this distribution is \emph{infinitely divisible}.

As our second main result, we obtain the following functional
non-central limit theorem (FNCLT) for generalized variations. The proof
of this result is carried out in Section~\ref{ssecnon-central} and
Section~\ref{ssectightness}.

\begin{thm}[(FNCLT)]\label{gvar-nclt}
Let $f$ be as above such that Assumption~\ref{Hermite-assumption} holds
and suppose that $H \in (1-\frac{1}{2\underline{k}},1 )^d$. Then
%nuo cia
\begin{equation}
\label{nclt-conv}
\overline{U}_f^{(n)} \mathop{\longrightarrow}_{n \rightarrow\infty}^{\mathbf{P}} \Lambda _{H,f}^{{1}/{2}}\widehat{Z}
 \qquad\mbox{in }D\bigl([0,1]^d\bigr),
\end{equation}
where $\widehat{Z}$ is a $d$-parameter Hermite sheet of order
$\underline{k}$ with Hurst parameter $\widetilde{H}$, given by \eqref
{H-shift}, and $\Lambda_{H,f}$ is given by \eqref{lambdadef}.
\end{thm}

\begin{rem}
(1) The convergence \eqref{nclt-conv} holds pointwise in $L^2(\Omega
,\mathscr{F},\mathbf{P})$, even when Assumption~\ref{Hermite-assumption}
does not hold; see Proposition~\ref{PNCLT} below.

(2) Unlike in Theorem~\ref{gvar-clt}, the non-central limit
$\widehat
{Z}$ is defined on the original probability space $(\Omega,\mathscr
{F},\mathbf{P})$. In particular, $\widehat{Z}$ is driven by the same white
noise $\mathscr{W}$ as $Z$.

(3) In the special case $\underline{k}=1$, the limit in \eqref
{nclt-conv} is Gaussian. In fact, then $\Lambda_{H,f} = a^2_1$ and
$\widehat{Z}=Z$.
\end{rem}

\begin{rem}
Our method of proving the convergence of finite-dimensional
distributions of $\overline{U}^{(n)}_f$, using chaotic expansions, is
particularly suitable for providing estimates on the \emph{speed of
convergence} (e.g., in the \emph{Wasserstein distance}) as is
done in \cite{NPR} following the original idea presented in \cite{NP},
which combines the \emph{Malliavin calculus} and \emph{Stein's method}.
In addition, the study of \emph{weighted} variations of the fBs is
still partially incomplete, especially with regards to functional
convergence (see \cite{Rev2009b}). To keep the length of this paper
within limits -- and since proving functional convergence of weighted
variations requires slightly different methods -- we have decided to
treat these two questions in a separate paper.
\end{rem}

\section{Finite-dimensional convergence}\label{fdd-convergence}

In this section, we begin the proofs of Theorems \ref{gvar-clt} and
\ref
{gvar-nclt}. To be more precise, we prove the finite-dimensional
statements corresponding to \eqref{clt-conv} and \eqref{nclt-conv}; see
Propositions \ref{fdd-conv-prop} and \ref{PNCLT}, respectively. As a
preparation, we study the correlation structure of the increments of
the fBs $Z$ and recall the chaotic expansion of functionals of $Z$.

\subsection{Correlation structure of increments}

In what follows, it will be convenient to use the shorthand
\begin{equation}
\label{increments}
Z^{(n)}_{i}:=\bigl\langle m(n)^H
\bigr\rangle Z \biggl( \biggl[\frac
{i-1}{m(n)},\frac{i}{m(n)} \biggr) \biggr),\qquad  1
\leq i \leq m(n),   n \in\mathbb{N}.
\end{equation}
For any $n \in\mathbb{N}$, the family $ \{Z^{(n)}_{i} \dvt 1 \leq i
\leq
m(n) \}$ is clearly centered and Gaussian. We will next derive its
correlation structure.

To describe the correlation structure of the rescaled increments \eqref
{increments}, let $\{B_{\check{H}}(t) \dvt t \in\mathbb{R}\}$ be an auxiliary
fractional Brownian motion with Hurst parameter $\check{H} \in(0,1)$.
Using the kernel \eqref{MVN}, we may represent it as
\begin{eqnarray*}
&& B_{\check{H}}(t) :=\int_\mathbb{R}G^{(1)}_{\check{H}}
(t,u) \,\mathrm{d}B(u), \qquad  t \in\mathbb{R},
\end{eqnarray*}
where $\{B(t) \dvt t \in\mathbb{R}\}$ is a standard Brownian motion.
Recall that
$B_{\check{H}}$ is $\check{H}$-\emph{self similar}, that is,
\begin{equation}
\label{ssi} \bigl\{B_{\check{H}}(at) \dvt t\in\mathbb{R} \bigr\} \stackrel{\mathscr{L}} {=} \bigl\{ a^{\check{H}}B_{\check{H}}(t) \dvt t\in\mathbb{R}
\bigr\} \qquad \mbox{for any }a > 0,
\end{equation}
and has \emph{stationary increments}, that is,
\begin{equation}
\label{si}
\bigl\{B_{\check{H}}\bigl([s,s+t)\bigr) \dvt t\in\mathbb{R} \bigr
\}\stackrel {\mathscr {L}}{=}\bigl\{ B_{\check{H}}\bigl([0,t)\bigr) \dvt t \in
\mathbb{R}\bigr\} \qquad \mbox{for any }s \in\mathbb{R}.
\end{equation}

The discrete parameter process
\begin{eqnarray*}
&& B_{\check{H}}\bigl([j,j+1\bigr)), \qquad  j \in\mathbb{Z},
\end{eqnarray*}
which is stationary by \eqref{si}, is called a \emph{fractional
Gaussian noise}.
Its correlation function can be expressed as
\begin{eqnarray*}
 r_{\check{H}}(j) &:=& \mathbf{E} \bigl[B_{\check{H}}\bigl([j,j+1
)\bigr)B_{\check
{H}}\bigl([0,1)\bigr) \bigr]\\
 &=& \frac{ |j+1|^{2\check{H}} -2|j|^{2\check
{H}}+|j-1|^{2\check{H}}}{2},  \qquad j \in
\mathbb{Z}.
\end{eqnarray*}
One can show, for example, using the mean value theorem, that there
exists a constant $C(\check{H}) > 0$
such that
\begin{equation}
\label{corr-decay}
\bigl|r_{\check{H}}(j)\bigr| \leq C(\check{H})|j|^{-2(1-\check{H})},\qquad  j
\in \mathbb{Z}.
\end{equation}
Thus, if $k > \frac{1}{2}$ and $\check{H} \in (0,1-\frac
{1}{2k}
)$, then
\begin{equation}
\label{corr-series}
\sum_{j \in\mathbb{Z}} \bigl|r_{\check{H}}(j)\bigr|^k
< \infty.
\end{equation}
If $\check{H} \in [1-\frac{1}{2k},1 )$, then the series
\eqref
{corr-series} diverges. In this case, it is still useful to have
estimates for the partial sums corresponding to \eqref{corr-series}.
Using \eqref{corr-decay}, one can prove that there exists a constant
$C'(\check{H},k)>0$ such that
\begin{equation}
\label{divergent-correlations}
\sum_{|j|<l} \bigl|r_{\check{H}}(j)\bigr|^k
\leq %
\cases{\displaystyle C'(\check{H},k) \log l, & \quad$\displaystyle\check{H} =
1-\frac{1}{2k}$,
\vspace*{3pt}\cr
\displaystyle C'(\check{H},k) l^{1-2k(1-\check{H})}, & \quad$\displaystyle\check{H} \in \biggl(1-
\frac
{1}{2k},1 \biggr)$.}
\end{equation}

We can now describe the correlations of the rescaled increments \eqref
{increments} using the correlation function of the fractional Gaussian
noise as follows.

\begin{lem}[(Correlation structure)]\label{corr-struct} For any $n \in
\mathbb{N}
$, and $1 \leq i^{(1)},  i^{(2)} \leq m(n)$,
\begin{eqnarray*}
&& \mathbf{E} \bigl[Z^{(n)}_{i^{(1)}}Z^{(n)}_{i^{(2)}}
\bigr] = \prod_{\nu=1}^d r_{H_\nu}
\bigl(i^{(1)}_\nu-i^{(2)}_\nu \bigr).
\end{eqnarray*}
\end{lem}

\begin{pf} Using first the linearity of Wiener integrals and then
the product structure \eqref{tensor-kernel} of the kernel $G_{H}^{(d)}$
and Remark~\ref{tensorcase}, we obtain
for any $s$,  $t\in[0,1]^d$ such that $s \leq t$,
\begin{eqnarray}
Z\bigl([s,t)\bigr) &=& \int G^{(d)}_{H}
\bigl([s,t),u\bigr) \mathscr{W}(\mathrm{d}u)
\nonumber
\\[-8pt]
\label{kernel-diff}
\\[-8pt]
\nonumber
&=& \int \prod
_{\nu=1}^d G^{(1)}_{H_\nu}
\bigl([s_\nu,t_\nu),u_\nu \bigr) \mathscr {W}(
\mathrm{d}u).
\end{eqnarray}
Thus, by Fubini's theorem,
\begin{eqnarray*}
&& \mathbf{E} \biggl[Z \biggl( \biggl[\frac{i^{(1)}-1}{m(n)},\frac
{i^{(1)}}{m(n)} \biggr) \biggr) Z \biggl(
\biggl[\frac{i^{(2)}-1}{m(n)},\frac
{i^{(2)}}{m(n)} \biggr) \biggr) \biggr]
\\
&&\quad= \prod_{\nu=1}^d
\int G^{(1)}_{H_\nu} \biggl( \biggl[\frac
{i^{(1)}_\nu
-1}{m_\nu(n)},
\frac{i^{(1)}_\nu}{m_\nu(n)} \biggr),v \biggr)G^{(1)}_{H_\nu
} \biggl( \biggl[
\frac{i^{(2)}_\nu-1}{m_\nu(n)},\frac{i^{(2)}_\nu
}{m_\nu
(n)} \biggr),v \biggr)\, \mathrm{d}v
\\
&&\quad = \prod_{\nu=1}^d \mathbf{E} \biggl[
B_{H_\nu} \biggl( \biggl[\frac
{i^{(1)}_\nu
-1}{m_\nu(n)},\frac{i^{(1)}_\nu}{m_\nu(n)} \biggr)
\biggr)B_{H_\nu
} \biggl( \biggl[\frac{i^{(2)}_\nu-1}{m_\nu(n)},\frac{i^{(2)}_\nu}{m_\nu
(n)} \biggr) \biggr)
\biggr].
\end{eqnarray*}
For any $\nu\in\{1,\ldots,d\}$, the fractional Brownian motion
$B_{H_\nu}$ is $H_\nu$-self similar and has stationary increments,
cf.
\eqref{ssi} and \eqref{si}, so we obtain
\begin{eqnarray*}
&& \mathbf{E} \biggl[ B_{H_\nu} \biggl( \biggl[\frac{i^{(1)}_\nu-1}{m_\nu
(n)},\frac
{i^{(1)}_\nu}{m_\nu(n)}
\biggr) \biggr)B_{H_\nu} \biggl( \biggl[\frac
{i^{(2)}_\nu-1}{m_\nu(n)},\frac{i^{(2)}_\nu}{m_\nu(n)}
\biggr) \biggr) \biggr]
= \frac{r_{H_\nu} (i^{(1)}_\nu-i^{(2)}_\nu )}{m_\nu
(n)^{2H_\nu}},
\end{eqnarray*}
from which the assertion follows.
\end{pf}

\subsection{Multiple Wiener integrals and central limit theorem}\label
{sseccltproof}

The proofs of Theorems \ref{gvar-clt} and \ref{gvar-nclt} rely on
particular representations of generalized variations in terms of \emph
{multiple Wiener integrals} with respect to the underlying white noise
$\mathscr{W}$. We will now briefly review the theory of multiple Wiener
integrals and how these integrals can be used to prove central limit
theorems. As an application, we take the first step in the proof of
Theorem~\ref{gvar-clt} by establishing the convergence of
finite-dimensional laws.

In what follows, we write $\mathscr{H} :=L^2(\mathbb{R}^d)$. Recall that
$\mathscr{H}$ is a separable Hilbert space when we endow it with the
usual inner product. For any $k\in\mathbb{N}$, we denote by $\mathscr
{H}^{\otimes k}$ the $k$-fold tensor power of $\mathscr{H}$ and by
$\mathscr{H}^{\odot k} \subset\mathscr{H}^{\otimes k}$ the
symmetrization of $\mathscr{H}^{\otimes k}$. Note that we can make the
identification
$\mathscr{H}^{\otimes k} \cong L^2(\mathbb{R}^{kd})$.
For any $h \in\mathscr{H}^{\odot k}$, we may define the $k$-fold
multiple Wiener integral $I^{\mathscr{W}}_k(h)$ of $h$ with respect to
$\mathscr{W}$. This is done, using Hermite polynomials, by setting for
any $\kappa\in\mathbb{N}$, any orthonormal $h_1,\ldots,h_\kappa\in
\mathscr
{H}$, and for any $k_1,\ldots,k_\kappa\in\mathbb{N}$ such that
$k_1+\cdots
+k_\kappa= k$,
\begin{equation}
\label{MW-def}
I^{\mathscr{W}}_k \Biggl(\bigodot_{j = 1}^\kappa
h_j^{\otimes
k_j} \Biggr) :=k! \prod_{j=1}^\kappa
P_{k_j} \biggl(\int h_j(u)\mathscr {W}(\mathrm{d} u)
\biggr),
\end{equation}
where $\odot$ denotes symmetrization of the tensor product, and
extended to general integrands $h \in\mathscr{H}^{\odot k}$ using a
density argument. It is worth stressing that the multiple Wiener
integral is linear with respect to the integrand and has zero
expectation. Moreover, by \eqref{MW-def}, for $h \in\mathscr{H}$ one has
\begin{equation}
\label{Wiener-i}
I^{\mathscr{W}}_1(h) = \int h(u) \mathscr{W}(
\mathrm{d}u),
\end{equation}
and if $\| h \|_{\mathscr{H}} =1$, then for any $k \in\mathbb{N}$ it holds
that $h^{\otimes k} \in\mathscr{H}^{\odot k}$\vspace*{-2pt} and
\begin{equation}
\label{Hermite-Wiener} P_k \bigl(I^{\mathscr{W}}_1(h) \bigr) =
I^{\mathscr{W}}_k\bigl(h^{\otimes k}\bigr).
\end{equation}
Multiple Wiener integrals have the following isometry and orthogonality
properties: for any $k_1$,  $k_2 \in\mathbb{N}$, $h_1\in\mathscr
{H}^{\odot
k_1}$, and $h_2\in\mathscr{H}^{\odot k_2}$,
\begin{equation}
\label{isometry}
\mathbf{E} \bigl[I^{\mathscr{W}}_{k_1}(h_1)I^{\mathscr
{W}}_{k_2}(h_2)
\bigr] = %
\cases{ k_1 ! \langle h_1,h_2
\rangle_{\mathscr{H}^{\otimes k_1}}, & \quad $k_1=k_2$,
\vspace*{3pt}\cr
0,& \quad $k_1\neq k_2$.}
\end{equation}

Recall that any random variable $Y \in L^2(\Omega,\mathscr{F},\mathbf{P})$
has a \emph{chaotic expansion} in terms of kernels $F^Y_k \in\mathscr
{H}^{\odot k}$, $k \in\mathbb{N}$,\vspace*{-3pt} as
\begin{equation}
\label{chaotic}
Y = \mathbf{E}[Y]+\sum_{k=1}^\infty
I^{\mathscr{W}}_k \bigl(F^Y_k \bigr),
\end{equation}
where the series converges in $L^2(\Omega,\mathscr{F},\mathbf{P})$ (see,
e.g., \cite{Kal2002}, Theorem~13.26). Since the apperance of the seminal
paper of Nualart and Peccati \cite{NP2005}, the convergence of random
variables admitting expansions of the form \eqref{chaotic} to a
Gaussian law has been well understood, based on convenient
characterizations using the properties of the kernels.
To describe the key result, recall that for any $k_1$,  $k_2$,  $r
\in\mathbb{N}$ such the $r < \min\{k_1,k_2\}$, the $r$th \emph
{contraction} of
$h_1 \in\mathscr{H}^{\otimes k_1}$ and $h_2 \in\mathscr{H}^{\otimes
k_2}$ is defined\vspace*{-2pt} as
\begin{eqnarray*}
&& (h_1 \otimes_r h_2) \bigl(t^{(1)},
\ldots,t^{(k_1+k_2-2r)} \bigr)
\\[-2pt]
&&\quad:= \bigl\langle h_1 \bigl(t^{(1)},\ldots,
t^{(k_1-r)},   \cdot  \bigr), h_2 \bigl(  \cdot
 ,t^{(k_1-r+1)},\ldots, t^{(k_1+k_2-2r)} \bigr) \bigr\rangle_{\mathscr{H}^{\otimes r}}
\end{eqnarray*}
for any $t^{(1)},\ldots,t^{(k_1+k_2-2r)} \in\mathbb{R}^d$. The following
multivariate central limit theorem for chaotic expansions appears in
\cite{BCP2009}, Theorem~5, where it is proven using the results in
\cite
{PT2005}.

\begin{lem}[(CLT for chaotic expansions)]\label{chaosclt} Let $\kappa
\in
\mathbb{N}$ and suppose that for any $n \in\mathbb{N}$, we are given
random variables
$Y^{(n)}_1,\ldots,Y^{(n)}_\kappa\in L^2(\Omega, \mathscr{F},\mathbf{P})$
such that\vspace*{-2pt} for any $j \in\{1,\ldots,\kappa\}$,
\begin{eqnarray*}
&& Y^{(n)}_j = \sum_{k=1}^{\infty}
I^{\mathscr{W}}_k \bigl(F^{(n)}_{k}(j,\cdot)
\bigr),
\end{eqnarray*}
where $F^{(n)}_{k}(j,\cdot) \in\mathscr{H}^{\odot k}$, $k \in
\mathbb{N}$. Let us
assume that the following conditions\vspace*{-2pt} hold:
\begin{longlist}[(a)]
\item[(a)] For\vspace*{-3pt} any $j \in\{1,\ldots,\kappa\}$,
\begin{eqnarray*}
&& \limsup_{n\rightarrow\infty}\sum_{k=K}^{\infty}
k! \bigl\| F^{(n)}_{k}(j,\cdot) \bigr\|^2_{\mathscr{H}^{\otimes k}}
\mathop{\longrightarrow}_{K \rightarrow\infty} 0.
\end{eqnarray*}
\item[(b)] There exists a sequence $\Sigma,\Sigma_1,\Sigma
_2,\ldots$ of positive semidefinite $d \times d$-matrices such that for
any $(j_1,j_2) \in\{1,\ldots,\kappa\}^2$ and\vspace*{-2pt} $k \in\mathbb{N}$,
\begin{eqnarray*}
&& k! \bigl\langle F^{(n)}_{k}(j_1,
\cdot),F^{(n)}_{k}(j_2,\cdot) \bigr\rangle
_{\mathscr{H}^{\otimes k}} \mathop{\longrightarrow}_{n \rightarrow\infty} \Sigma_k(j_1,j_2),
\end{eqnarray*}
and that $\sum_{k=1}^\infty\Sigma_k = \Sigma$.
\item[(c)] For any $j \in\{1,\ldots,\kappa\}$, $k\geq2$, and
$r \in\{1,\ldots,k-1\}$,
\begin{eqnarray*}
&& \bigl\| F^{(n)}_{k}(j,\cdot) \otimes_r
F^{(n)}_{k}(j,\cdot) \bigr\| ^2_{\mathscr{H}^{\otimes2(k-r)}}
\mathop{\longrightarrow}_{n \rightarrow\infty} 0.
\end{eqnarray*}
\end{longlist}
Then we have
\begin{eqnarray*}
&& \bigl(Y^{(n)}_1,\ldots,Y^{(n)}_\kappa
\bigr) \mathop{\longrightarrow}_{n \rightarrow \infty}^{\mathscr{L}} N_\kappa(0,
\Sigma),
\end{eqnarray*}
where $N_\kappa(0,  \Sigma)$ stands for the $\kappa$-dimensional
Gaussian law with mean $0$ and covariance matrix~$\Sigma$.
\end{lem}

We apply now Lemma~\ref{chaosclt} to establish the following
finite-dimensional version of Theorem~\ref{gvar-clt}.

\begin{prop}[(CLT for finite-dimensional laws)]\label{fdd-conv-prop}
Suppose that $H \in(0,1)^d \setminus (1-\frac{1}{2\underline
{k}},1 )^d$.
Let $\kappa\in\mathbb{N}$ and $ (t^{(1)},\ldots,t^{(\kappa
)} ) \in
([0,1]^d )^\kappa$. Then
\begin{equation}
\label{fdd-conv}
\bigl(Z\bigl(t^{(1)}\bigr),\ldots,Z\bigl(t^{(\kappa)}
\bigr),\overline {U}^{(n)}_f\bigl(t^{(1)}\bigr),
\ldots,\overline{U}^{(n)}_f\bigl(t^{(\kappa)}\bigr)
\bigr) \mathop{\longrightarrow}_{n\rightarrow\infty}^{\mathscr{L}} N_{2\kappa} \left(0,
\left[\matrix{ \Xi& 0
\cr
0 & \Sigma
} \right]\right),
\end{equation}
where $\Xi$ is the covariance matrix of the random vector $
(Z(t^{(1)}),\ldots,Z(t^{(\kappa)}) )$ and
\begin{eqnarray*}
&& \Sigma(j_1,j_2) :=\Lambda_{H,f}
R^{(d)}_{\widetilde
{H}}\bigl(t^{(j_1)},t^{(j_2)}\bigr),
\qquad (j_1,j_2) \in\{1,\ldots,\kappa\}^2.
\end{eqnarray*}
\end{prop}

\begin{rem}
In the case $H \in (0,1-\frac{1}{2\underline{k}} ]^d$, the
convergence
\begin{eqnarray*}
\bigl(\overline{U}^{(n)}_f\bigl(t^{(1)}\bigr),
\ldots,\overline {U}^{(n)}_f\bigl(t^{(\kappa)}\bigr)
\bigr) \mathop{\longrightarrow}_{n\rightarrow\infty }^{\mathscr{L}} N_{\kappa}(\Sigma)
\end{eqnarray*}
follows from the classical results of Breuer and Major \cite{BM}.
\end{rem}

\begin{pf*}{Proof of Proposition~\ref{fdd-conv-prop}}
By \eqref{Wiener-i}, we have $Z(t) = I^{\mathscr{W}}_1
(G^{(d)}_H(t,\cdot)  )$ for any $t \in[0,1]^d$. In particular, by
\eqref{kernel-diff} and linearity, we find that for any $n \in\mathbb
{N}$ and
$1 \leq i \leq m(n)$,
\begin{eqnarray*}
&& Z^{(n)}_i = I^{\mathscr{W}}_1
\bigl(h^{(n)}_i \bigr),
\end{eqnarray*}
where
\begin{equation}
\label{hgdef}
h^{(n)}_i :=\bigl\langle m(n)^H
\bigr\rangle g^{(n)}_i,\qquad   g^{(n)}_i
:=G^{(d)}_H \biggl( \biggl[\frac{i-1}{m(n)},\frac{i}{m(n)}
\biggr),  \cdot  \biggr),
\end{equation}
satisfying $ \|h^{(n)}_i \|_\mathscr{H}=1$, due to the relation
\eqref{isometry} and Lemma~\ref{corr-struct}.
The expansion \eqref{f-Hermite} and the connection of Hermite
polynomials and multiple Wiener integrals \eqref{Hermite-Wiener} allows
then us to write
\begin{equation}
\label{U-chaotic}
\overline{U}^{(n)}_f(t) = \sum
_{k = \underline{k}}^\infty I^{\mathscr
{W}}_k
\bigl(F^{(n)}_k(t,\cdot) \bigr),  \qquad t
\in[0,1]^d,   n \in \mathbb{N},
\end{equation}
where
\begin{equation}
\label{F-def} F^{(n)}_k(t,\cdot) :=\frac{a_k}{\langle c(n) \rangle^{1/2}} \sum
_{1 \leq i \leq\lfloor m(n)t \rfloor} \bigl(h^{(n)}_i
\bigr)^{\otimes k},  \qquad k \geq\underline{k}.
\end{equation}

For the remainder of the proof, let $s$, $t \in\{t^{(1)},\ldots
,t^{(\kappa)}\}$.
Let us first look into condition (a) of Lemma~\ref{chaosclt}.
By Lemma~\ref{corr-struct} and the relation \eqref{isometry}, we obtain
for any $n \in\mathbb{N}$ and $k \geq\underline{k}$,
\begin{eqnarray}
\nonumber
\bigl\langle F^{(n)}_k(s,
\cdot),F^{(n)}_k(t,\cdot) \bigr\rangle _{\mathscr
{H}^{\otimes k}} & =&
\frac{a^2_k}{\langle c(n) \rangle} \sum_{1 \leq
i^{(1)} \leq\lfloor m(n)s \rfloor} \sum
_{1 \leq i^{(2)} \leq\lfloor
m(n)t \rfloor} \bigl\langle \bigl(h^{(n)}_{i^{(1)}}
\bigr)^{\otimes k}, \bigl(h^{(n)}_{i^{(2)}}
\bigr)^{\otimes k} \bigr\rangle_{\mathscr
{H}^{\otimes
k}}
\\
\label{innerprod}
& =&  \frac{a^2_k}{\langle c(n)\rangle} \sum_{1 \leq i^{(1)} \leq
\lfloor
m(n)s \rfloor} \sum
_{1 \leq i^{(2)} \leq\lfloor m(n)t \rfloor} \bigl\langle h^{(n)}_{i^{(1)}},
h^{(n)}_{i^{(2)}} \bigr\rangle_{\mathscr
{H}}^k
\\
&=&  a^2_k \prod_{\nu= 1}^d
\frac{1}{c_\nu(n)} \sum_{j_1=1}^{\lfloor
m_\nu(n) s_\nu\rfloor} \sum
_{j_2=1}^{\lfloor m_\nu(n) t_\nu
\rfloor} r_{H_\nu}(j_1-j_2)^k.
\nonumber
\end{eqnarray}
Let $k_0 \in\mathbb{N}$ be large enough so that $H_\nu\in
(0,1-\frac
{1}{2k_0} )$ for any $\nu\in\{1,\ldots,d\}$. Then we have for any
$k \geq k_0$,
\begin{eqnarray*}
0 &\leq & \prod_{\nu= 1}^d
\frac{1}{c_\nu(n)} \sum_{j_1 =1}^{\lfloor
m_\nu
(n) s_\nu\rfloor} \sum
_{j_2 =1}^{\lfloor m_\nu(n) s_\nu\rfloor} r_{H_\nu}(j_1-j_2)^k
\\
&\leq & \prod_{\nu= 1}^d \frac{1}{c_\nu(n)}
\sum_{j_1=1}^{m_\nu(n)} \sum
_{j_2=1}^{m_\nu(n)} \bigl|r_{H_\nu}(j_1-j_2)\bigr|^{k_0}
\\
& \leq & \prod_{\nu= 1}^d \biggl(\sup
_{n \in\mathbb{N}} \frac{m_\nu
(n)}{c_\nu
(n)} \biggr) \sum
_{j \in\mathbb{Z}}\bigl|r_{H_\nu}(j)\bigr|^{k_0} <\infty,
\end{eqnarray*}
which follows from Remark~\ref{psidelta-comparison} and the elementary estimate
\begin{equation}
\label{sup-sum}
\sup_{1\leq j_1\leq\ell} \sum_{j_2=1}^\ell\bigl|r_{H_\nu}(j_1-j_2)\bigr|^q
\leq \sum_{|j|<\ell}\bigl|r_{H_\nu}(j)\bigr|^q,
 \qquad\ell\in\mathbb{N},   q \in\mathbb{R}_+.
\end{equation}
Thus, by \eqref{standard-sum}, we have for $K \geq k_0$,
\begin{eqnarray*}
&& 0 \leq\limsup_{n \rightarrow\infty}\sum_{k=K}^{\infty}
k! \bigl\| F^{(n)}_k(s) \bigr\|^2_{\mathscr{H}^{\otimes k}} \leq\sum
_{k=K}^{\infty} k! a^2_k
\prod_{\nu= 1}^d \biggl(\sup
_{n \in\mathbb{N}} \frac
{m_\nu(n)}{
c_\nu
(n)} \biggr) \sum
_{j \in\mathbb{Z}}\bigl|r_{H_\nu}(j)\bigr|^{k_0} \mathop{\longrightarrow}_{K \rightarrow\infty} 0,
\end{eqnarray*}
and the condition (a) is verified.

To check condition (b) of Lemma~\ref{chaosclt}, note that we
can write for any $\nu\in\{1,\ldots,d\}$, assuming without loss of
generality that $t_\nu\geq s_\nu$,
\begin{eqnarray}
\nonumber
&& \frac{1}{c_\nu(n)} \sum_{j_1=1}^{\lfloor m_\nu(n) s_\nu\rfloor}
\sum_{j_2=1}^{\lfloor m_\nu(n) t_\nu\rfloor} r_{H_\nu}(j_1-j_2)^k
\\
\label{fixed-nu}
&&\quad = \frac{1}{2} \Biggl(\frac{1}{c_\nu(n)} \sum
_{j_1=1}^{\lfloor
m_\nu(n)
s_\nu\rfloor} \sum
_{j_2=1}^{\lfloor m_\nu(n) s_\nu\rfloor} r_{H_\nu
}(j_1-j_2)^k
+\frac{1}{c_\nu(n)} \sum_{j_1=1}^{\lfloor m_\nu(n)
t_\nu
\rfloor} \sum
_{j_2=1}^{\lfloor m_\nu(n) t_\nu\rfloor} r_{H_\nu
}(j_1-j_2)^k\qquad\quad
\\
\nonumber
&&\hspace*{14pt}\qquad  {}- \frac{1}{c_\nu(n)} \sum_{j_1=1}^{\lfloor m_\nu(n)
t_\nu\rfloor-\lfloor m_\nu(n) s_\nu\rfloor}
\sum_{j_2=1}^{\lfloor
m_\nu(n) t_\nu\rfloor-\lfloor m_\nu(n) s_\nu\rfloor} r_{H_\nu
}(j_1-j_2)^k
\Biggr).
\end{eqnarray}
We will now compute the limit of \eqref{fixed-nu} separately in the
following three possible cases:
\begin{longlist}[(iii)]
\item[(i)] $H_\nu\in (1-\frac{1}{2\underline{k}},1 )$,
\item[(ii)] $H_\nu= 1-\frac{1}{2\underline{k}}$,
\item[(iii)] $H_\nu\in (0,1-\frac{1}{2\underline{k}} )$.
\end{longlist}
In the case (i), we obtain, by Lemma A.1 of \cite{RST2012},
\begin{eqnarray*}
&& \frac{1}{c_\nu(n)} \sum_{j_1=1}^{\lfloor m_\nu(n) s_\nu\rfloor} \sum
_{j_2=1}^{\lfloor m_\nu(n) s_\nu\rfloor} r_{H_\nu
}(j_1-j_2)^{\underline
{k}}
\\
&&\quad= \biggl(\frac{\lfloor m_\nu(n) s_\nu\rfloor}{m_\nu(n)} \biggr)^{2-2\underline{k}(1-H_\nu)}\bigl\lfloor
m_\nu(n) s_\nu\bigr\rfloor ^{-2+2\underline{k}(1-H_\nu)} \sum
_{j_1=1}^{\lfloor m_\nu(n) s_\nu
\rfloor} \sum_{j_2=1}^{\lfloor m_\nu(n) s_\nu\rfloor}
r_{H_\nu
}(j_1-j_2)^{\underline{k}}
\\
&&\quad\mathop{\longrightarrow}_{n \rightarrow\infty} \kappa(H_\nu,\underline
{k})s_\nu ^{2-2\underline{k}(1-H_\nu)} = \kappa(H_\nu,
\underline{k})s_\nu ^{\widetilde{H}_\nu},
\end{eqnarray*}
where $\kappa(H_\nu,\underline{k})$ is given by \eqref{bdef}.
(In fact, Lemma A.1 of \cite{RST2012} requires that $\underline
{k}\geq
2$, but it is straightforward to check that the limits stated therein
are valid also when $\underline{k}=1$.)
With $k > \underline{k}$ we may choose $\varepsilon>0$ so that
$\underline{k}+\varepsilon< \min (\frac{1}{2(1-H_\nu)},k
)$, whence
\begin{eqnarray}
\Biggl|\frac{1}{c_\nu(n)} \sum
_{j_1=1}^{\lfloor m_\nu(n) s_\nu
\rfloor} \sum_{j_2=1}^{\lfloor m_\nu(n) s_\nu\rfloor}
r_{H_\nu
}(j_1-j_2)^k \Biggr| \!\!&& \leq
\frac{1}{c_\nu(n)} \sum_{j_1=1}^{m_\nu(n)} \sum
_{j_2=1}^{m_\nu
(n)} \bigl|r_{H_\nu}(j_1-j_2)\bigr|^{\underline{k}+\varepsilon}
\nonumber
\\
\label{epsilon-power}
&& \leq \frac{1}{m_\nu(n)^{1-2\underline{k}(1-H_\nu)}}\sum_{|j|<m_\nu
(n)}
\bigl|r_{H_\nu}(j)\bigr|^{\underline{k}+\varepsilon}\qquad\quad \\
\nonumber
&&  \displaystyle\mathop{\longrightarrow}_{n \rightarrow\infty}  0
\end{eqnarray}
by the estimate \eqref{divergent-correlations}.
Treating the other summands on the right-hand side of \eqref{fixed-nu}
similarly, we arrive at
\begin{eqnarray*}
&&\lim_{n \rightarrow\infty}\frac{1}{c_\nu(n)} \sum
_{j_1=1}^{\lfloor
m_\nu(n) s_\nu\rfloor} \sum_{j_2=1}^{\lfloor m_\nu(n) t_\nu
\rfloor}
r_{H_\nu}(j_1-j_2)^k
\\
&&\quad=
\cases{ \displaystyle\frac{\kappa(H_\nu,\underline{k})}{2} \bigl(s_\nu^{\widetilde
{H}_\nu
}+t_\nu^{\widetilde{H}_\nu}
- (t_\nu-s_\nu)^{\widetilde{H}_\nu
} \bigr)=
\kappa(H_\nu,\underline{k})R^{(1)}_{\widetilde{H}_\nu}(s_\nu
,t_\nu), &\quad $k= \underline{k}$,\vspace*{3pt}
\cr
0, & \quad $k > \underline{k}$.}
\end{eqnarray*}
In the case (ii), rearranging and applying Lemma A.1 of \cite
{RST2012} yields
\begin{eqnarray*}
&& \frac{1}{c_\nu(n)} \sum_{j_1=1}^{\lfloor m_\nu(n) s_\nu\rfloor} \sum
_{j_2=1}^{\lfloor m_\nu(n) s_\nu\rfloor} r_{H_\nu
}(j_1-j_2)^{\underline
{k}}
\\
&&\quad= \biggl(1+\frac{\log ({\lfloor m_\nu(n) s_\nu\rfloor
}/{m_\nu
(n)}  )}{\log (m_\nu(n) )} \biggr) \frac{\lfloor m_\nu(n)
s_\nu\rfloor/m_\nu(n)}{\lfloor m_\nu(n) s_\nu\rfloor\log(\lfloor
m_\nu(n) s_\nu\rfloor)} \sum
_{j_1=1}^{\lfloor m_\nu(n) s_\nu
\rfloor} \sum_{j_2=1}^{\lfloor m_\nu(n) s_\nu\rfloor}
r_{H_\nu
}(j_1-j_2)^{\underline{k}}
\\
&&\quad\mathop{\longrightarrow}_{n \rightarrow\infty} \iota(\underline{k})s_\nu,
\end{eqnarray*}
where $\iota(\underline{k})$ is given by \eqref{bdef}.
When $k>\underline{k}$, we have $H_\nu\in (0,1-\frac{1}{2k} )$
and, consequently,
\begin{eqnarray*}
&& \Biggl|\frac{1}{c_\nu(n)} \sum_{j_1=1}^{\lfloor m_\nu(n) s_\nu
\rfloor} \sum
_{j_2=1}^{\lfloor m_\nu(n) s_\nu\rfloor} r_{H_\nu
}(j_1-j_2)^{k}
\Biggr|\leq\frac{1}{\log (m_\nu(n) )} \sum_{j \in\mathbb{Z}
}
\bigl|r_{H_\nu}(j)\bigr|^{k} \mathop{\longrightarrow}_{n \rightarrow\infty} 0.
\end{eqnarray*}
Again, a similar treatment of the other summands on right-hand side of
\eqref{fixed-nu} establishes that
\begin{eqnarray*}
&& \lim_{n \rightarrow\infty} \frac{1}{c_\nu(n)} \sum
_{j_1=1}^{\lfloor
m_\nu(n) s_\nu\rfloor} \sum_{j_2=1}^{\lfloor m_\nu(n) t_\nu
\rfloor}
r_{H_\nu}(j_1-j_2)^k
\\
&&\quad=
\cases{\displaystyle \frac{\iota(\underline{k})}{2} \bigl(s_\nu+
t_\nu-(t_\nu-s_\nu ) \bigr) = \iota(
\underline{k})R^{(1)}_{\widetilde{H}_\nu}(s_\nu,t_\nu),
&\quad $k= \underline{k}$,
\vspace*{3pt}\cr
0, & \quad $k > \underline{k}$.}
\end{eqnarray*}
Finally, in the case (iii), we deduce in a straightforward
manner that for any $k \geq\underline{k}$,
\begin{eqnarray*}
\lim_{n \rightarrow\infty} \frac{1}{c_\nu(n)} \sum
_{j_1=1}^{\lfloor
m_\nu(n) s_\nu\rfloor} \sum
_{j_2=1}^{\lfloor m_\nu(n) t_\nu
\rfloor} r_{H_\nu}(j_1-j_2)^k
& = & \frac{1}{2}\sum_{j \in\mathbb{Z}} r_{H_\nu}(j)^k
\bigl(s_\nu+ t_\nu-(t_\nu-s_\nu)
\bigr)
\\
& =& \sum_{j \in\mathbb{Z}} r_{H_\nu}(j)^k
R^{(1)}_{\widetilde
{H}_\nu
}(s_\nu ,t_\nu)
\end{eqnarray*}
using Lemma A.1 of \cite{RST2012}.

Returning to the expression \eqref{innerprod}, we have shown that for
any $k \geq\underline{k}$,
\begin{equation}
\label{cov-conv}
k! \bigl\langle F^{(n)}_k(s),F^{(n)}_k(t)
\bigr\rangle_{\mathscr
{H}^{\otimes k}}\mathop{\longrightarrow}_{n \rightarrow\infty} k!
a_k^2 \bigl\langle b^{(k)}\bigr\rangle
R^{(d)}_{\widetilde{H}}(s,t).
\end{equation}
When $\underline{k}=1$, we need to check, additionally, that the
covariance matrix appearing in the limit~\eqref{fdd-conv} is
block-diagonal. To this end, note that it follows from the assumption
$H \in(0,1)^d \setminus (\frac{1}{2},1 )^d$, that $b_\nu^{(1)}
= 0$ for some $\nu\in\{1,\ldots,d\}$, which in turn implies that
\begin{eqnarray*}
&& \bigl\|F^{(n)}_1(s,\cdot) \bigr\|^2_{\mathscr{H}}
\mathop{\longrightarrow}_{n \rightarrow\infty} 0.
\end{eqnarray*}
By the Cauchy--Schwarz inequality, we have then
\begin{eqnarray*}
&& \bigl\langle F^{(n)}_1(s,\cdot),G^{(d)}_H(t,
\cdot) \bigr\rangle _{\mathscr
{H}} \mathop{\longrightarrow}_{n \rightarrow\infty} 0,
\end{eqnarray*}
which ensures block diagonality, and concludes the verification of
condition (b).

In order to check condition (c) of Lemma~\ref{chaosclt}, let
$k \geq\max(\underline{k},2)$ and $r \in\{1,\ldots,k-1\}$. Using the
bilinearity of contractions and inner products, we obtain
\begin{eqnarray*}
&& \bigl\| F^{(n)}_{k}(t,\cdot) \otimes_r
F^{(n)}_{k}(t,\cdot) \bigr\| ^2_{\mathscr{H}^{\otimes2(k-r)}}
\\
&&\quad = \frac{a^4_k}{\langle c(n)\rangle^2} \mathop{\sum_{1 \leq
i^{(j)}\leq\lfloor m(n)t\rfloor}}_{j \in\{1,2,3,4\}} \bigl\langle \bigl(h^{(n)}_{i^{(1)}}
\bigr)^{\otimes k} \otimes_r \bigl(h^{(n)}_{i^{(2)}}
\bigr)^{\otimes k}, \bigl(h^{(n)}_{i^{(3)}}
\bigr)^{\otimes k} \otimes_r \bigl(h^{(n)}_{i^{(4)}}
\bigr)^{\otimes k} \bigr\rangle_{\mathscr
{H}^{\otimes2(k-r)}}
\\
&&\quad = \frac{a^4_k}{\langle c(n)\rangle^2} \mathop{\sum_{1 \leq i^{(j)}\leq\lfloor m(n)t\rfloor}}_{j \in\{1,2,3,4\}} \bigl\langle
h^{(n)}_{i^{(1)}},h^{(n)}_{i^{(2)}} \bigr
\rangle_{\mathscr{H}}^r \bigl\langle h^{(n)}_{i^{(3)}},h^{(n)}_{i^{(4)}}
\bigr\rangle_{\mathscr{H}}^r \bigl\langle h^{(n)}_{i^{(1)}},h^{(n)}_{i^{(3)}}
\bigr\rangle_{\mathscr
{H}}^{k-r} \bigl\langle h^{(n)}_{i^{(2)}},h^{(n)}_{i^{(4)}}
\bigr\rangle _{\mathscr{H}}^{k-r}
\\
&&\quad = a^4_k \prod_{\nu=1}^d
\frac{1}{c_\nu(n)^2} \sum_{j_1,j_2,j_3,j_4=1}^{\lfloor m_\nu(n) t_\nu\rfloor}
r_{H_\nu
}(j_1-j_2)^r
r_{H_\nu}(j_3-j_4)^r
r_{H_\nu}(j_1-j_3)^{k-r}
r_{H_\nu
}(j_2-j_4)^{k-r}.
\end{eqnarray*}
Following the proof of Lemma~4.1 of \cite{NPP2011}, we apply the bound
\begin{eqnarray*}
&& \bigl|r_{H_\nu}(j_1)\bigr|^r\bigl|r_{H_\nu}(j_2)\bigr|^{k-r}
\leq\bigl|r_{H_\nu}(j_1)\bigr|^k + \bigl|r_{H_\nu}(j_2)\bigr|^{k},
\qquad  j_1,  j_2 \in\mathbb{Z},
\end{eqnarray*}
which is a consequence of Young's inequality, and use repeatedly \eqref
{sup-sum} to deduce that
\begin{equation}
\label{J-product}
\bigl\| F^{(n)}_{k}(t,\cdot) \otimes_r
F^{(n)}_{k}(t,\cdot) \bigr\| ^2_{\mathscr{H}^{\otimes2(k-r)}}
\leq16^d a^4_k \prod
_{\nu=1}^d \frac{m_\nu(n)\phi_\nu(n)}{c_\nu
(n)^2} ,
\end{equation}
where
\begin{eqnarray*}
&& \phi_\nu(n) :=\sum_{|j_1|<m_\nu(n)} \bigl|r_{H_\nu}(j_1)\bigr|^r
\sum_{|j_2|<m_\nu(n)}\bigl|r_{H_\nu}(j_2)\bigr|^{k-r}
\sum_{|j_3|<m_\nu
(n)}\bigl|r_{H_\nu}(j_3)\bigr|^k.
\end{eqnarray*}

We need to analyze the asymptotic behaviour of $\phi_\nu(n)$ as $n
\rightarrow\infty$. This can be accomplished by considering separately
the three possible cases:
\begin{longlist}[(iii$'$)]
\item[(i$'$)] $H_\nu\in (1-\frac{1}{2k},1 )$,
\item[(ii$'$)] $H_\nu= 1-\frac{1}{2k}$,
\item[(iii$'$)] $H_\nu\in (0,1-\frac{1}{2k} )$.
\end{longlist}
In the case (i$'$) we have clearly $H_\nu\in (1-\frac
{1}{2(k-r)},1 ) \cap (1-\frac{1}{2r},1 )$, and by the
estimate \eqref{divergent-correlations}, it follows that
\begin{eqnarray*}
&& \phi_\nu(n) \leq C''(H_\nu,k,r)
m_\nu(n)^{3-4k(1-H_\nu)},
\end{eqnarray*}
where $C''(H_\nu,k,r) :=C'(H_\nu,r)C'(H_\nu,k-r)C'(H_\nu,k)$.
Since $H_\nu\in (1-\frac{1}{2k},1 ) \subset (1-\frac
{1}{2\underline{k}},1 )$, we obtain
\begin{eqnarray*}
&& \limsup_{n \rightarrow\infty}\frac{m_\nu(n)\phi_\nu(n)}{c_\nu(n)^2} \leq\limsup
_{n \rightarrow\infty}\frac{C''(H_\nu,k,r)}{m_\nu
(n)^{4(k-\underline{k})(1-H_\nu)}} < \infty.
\end{eqnarray*}
Let us then consider to the case (ii$'$). We have still $H_\nu
\in (1-\frac{1}{2(k-r)},1 ) \cap (1-\frac{1}{2r},1
)$, so
by \eqref{divergent-correlations} we find that
\begin{eqnarray*}
&& \phi_\nu(n) \leq C''(H_\nu,k,r)
m_\nu(n)^{2-2k(1-H_\nu)}\log \bigl( m_\nu (n) \bigr) =
C''(H_\nu,k,r) m_\nu(n)\log
\bigl( m_\nu(n) \bigr).
\end{eqnarray*}
Necessarily $H_\nu\in [1-\frac{1}{2\underline{k}},1 )$, whence
there is an index $n_0 \in\mathbb{N}$ such that $c_\nu(n) \geq
m_\nu(n) \log
 (m_\nu(n) )$ for all $n\geq n_0$. We deduce then that
\begin{eqnarray*}
&& \lim_{n \rightarrow\infty}\frac{m_\nu(n)\phi_\nu(n)}{c_\nu
(n)^2} \leq \lim_{n \rightarrow\infty}
\frac{C''(H_\nu,k,r)}{\log (m_\nu
(n)
)} =0.
\end{eqnarray*}
In the remaining case (iii$'$) we have $\sum_{j\in\mathbb{Z}}
|r_{H_\nu
}(j)|^k <\infty$. Since there is $n_0 \in\mathbb{N}$ such that
$c_\nu(n)
\geq
m_\nu(n)$ for all $n \geq n_0$, we find that
\begin{eqnarray*}
&& \lim_{n \rightarrow\infty}\frac{m_\nu(n)\phi_\nu(n)}{c_\nu
(n)^2}
\\
&&\quad\leq\lim_{n \rightarrow\infty}\frac{1}{m_\nu(n)}\sum
_{|j_1|<m_\nu
(n)}\bigl|r_{H_\nu}(j_1)\bigr|^r
\sum_{|j_2|<m_\nu(n)}\bigl|r_{H_\nu
}(j_2)\bigr|^{k-r}
\sum_{j_3\in\mathbb{Z}}\bigl|r_{H_\nu}(j_3)\bigr|^k
= 0
\end{eqnarray*}
by Lemma~2.2 of \cite{NPP2011}.

Finally, let us return to the upper bound \eqref{J-product}. The
crucial observation is that the assumption $H \in(0,1)^d \setminus
(1-\frac{1}{2\underline{k}},1 )^d$ implies that there is at least
one coordinate $\nu\in\{1,\ldots,d\}$ that falls within case (ii$'$) or (iii$'$). Thus,
\begin{eqnarray*}
&& \bigl\| F^{(n)}_{k}(t,\cdot) \otimes_r
F^{(n)}_{k}(t,\cdot) \bigr\| ^2_{\mathscr{H}^{\otimes2(k-r)}}
\mathop{\longrightarrow}_{n \rightarrow\infty} 0,
\end{eqnarray*}
concluding the verification of the condition (c), and the
convergence \eqref{fdd-conv} follows.
\end{pf*}

\subsection{Convergence to the Hermite sheet}\label{ssecnon-central}

We prove next a pointwise version of Theorem~\ref{gvar-nclt} in
$L^2(\Omega)$. The argument is based mainly on the chaotic expansion
\eqref{U-chaotic} and the isometry property \eqref{isometry} of
multiple Wiener integrals. However, compared to the proof of
Proposition~\ref{fdd-conv-prop}, we need to analyze the asymptotic
behaviour of the associated kernels more carefully.

\begin{prop}[(Pointwise NCLT)]\label{PNCLT}
Suppose that $H \in (1-\frac{1}{2\underline{k}},1 )^d$. Then,
for any $t \in[0,1]^d$,
\begin{equation}
\label{pointwise-L2}
\overline{U}^{(n)}_f(t) \mathop{\longrightarrow}_{n \rightarrow\infty }^{L^2(\Omega)} \Lambda_{H,f}^{{1}/{2}}
\widehat{Z}(t),
\end{equation}
where $\widehat{Z}$ is the Hermite sheet appearing in Theorem~\ref{gvar-nclt}.
\end{prop}

\begin{pf} Fix $t \in[0,1]^d$. By the chaotic expansion \eqref
{U-chaotic}, we have for any $n \in\mathbb{N}$,
\begin{eqnarray*}
&& \overline{U}^{(n)}_f(t) = I^{\mathscr{W}}_{\underline{k}}
\bigl(F^{(n)}_{\underline{k}}(t,\cdot) \bigr) + \sum
_{k = \underline
{k}+1}^\infty I^{\mathscr{W}}_k
\bigl(F^{(n)}_k(t,\cdot) \bigr).
\end{eqnarray*}
Using the property \eqref{isometry} and Parseval's identity, we find that
\begin{eqnarray*}
&& \mathbf{E} \Biggl[ \Biggl(\sum_{k = \underline{k}+1}^\infty
I^{\mathscr
{W}}_k \bigl(F^{(n)}_k(t,\cdot)
\bigr) \Biggr)^2 \Biggr] = \sum_{k = \underline
{k}+1}^\infty
k! \bigl\|F^{(n)}_k(t,\cdot) \bigr\|^2_{\mathscr
{H}^{\otimes k}}.
\end{eqnarray*}
Since $H \in (1-\frac{1}{2\underline{k}},1 )^d$, we may choose
$\varepsilon\in(0,1]$ so that $H \in (1-\frac{1}{2(\underline
{k}+\varepsilon)},1 )^d$. Combining \eqref{innerprod} and \eqref
{epsilon-power}, we find that
\begin{eqnarray*}
&& \sum_{k = \underline{k}+1}^\infty k! \bigl\|F^{(n)}_k(t,
\cdot) \bigr\| ^2_{\mathscr{H}^{\otimes k}} \leq\sum_{k = \underline{k}+1}^\infty
k! a^2_k \prod_{\nu=1}^d
\frac{1}{m_\nu(n)^{1-2\underline{k}(1-H_\nu
)}}\sum_{|j|<m_\nu(n)} \bigl|r_{H_\nu}(j)\bigr|^{\underline{k}+\varepsilon}
\mathop{\longrightarrow}_{n \rightarrow\infty} 0,
\end{eqnarray*}
where convergence to zero is a consequence of the bound \eqref
{divergent-correlations}. Thus, it remains to show that
\begin{eqnarray*}
&& I^{\mathscr{W}}_{\underline{k}} \bigl(F^{(n)}_{\underline{k}}(t,\cdot)
\bigr) \mathop{\longrightarrow}_{n \rightarrow\infty}^{L^2(\Omega)}I^{\mathscr
{W}}_{\underline{k}}
\bigl(\Lambda_{H,f}^{{1}/{2}} \widehat {G}^{(\underline{k})}_{\widetilde{H}}(t,
\cdot) \bigr) = \Lambda _{H,f}^{{1}/{2}} \widehat{Z}(t),
\end{eqnarray*}
which follows by \eqref{isometry}, if we can show that
\begin{equation}
\label{kernel-convergence}
F^{(n)}_{\underline{k}}(t,\cdot) \mathop{\longrightarrow}_{n \rightarrow\infty }^{\mathscr{H}^{\otimes\underline{k}}} \Lambda_{H,f}^{{1}/{2}}
\widehat{G}^{(\underline{k})}_{\widetilde{H}}(t,\cdot).
\end{equation}
In the special case $\underline{k}=1$, the convergence \eqref
{pointwise-L2} follows already. Namely,
\begin{eqnarray*}
&& I^{\mathscr{W}}_1 \bigl(F^{(n)}_1(t,\cdot)
\bigr) = a_1 Z \biggl( \frac
{\lfloor m(n) t \rfloor}{m(n)} \biggr) \mathop{\longrightarrow}_{n\rightarrow \infty }^{L^2(\Omega)} a_1 Z(t) =
\Lambda^{{1}/{2}}_{H,f} Z(t) = \Lambda ^{{1}/{2}}_{H,f}
\widehat{Z}(t)
\end{eqnarray*}
by the $L^2$-continuity of $Z$. Thus, we can assume that $\underline
{k}\geq2$ from now on.

We will prove the convergence \eqref{kernel-convergence} in two steps.
First, we show that $ (F^{(n)}_{\underline{k}}(t,\cdot) )_{n
\in
\mathbb{N}}$ is a Cauchy sequence in $\mathscr{H}^{\otimes\underline{k}}$.
Later, we characterize the limit. Let $n_1$,  $n_2 \in\mathbb{N}$
and consider
\begin{eqnarray}
&&\bigl\|F^{(n_1)}_{\underline{k}}(t,\cdot)-F^{(n_2)}_{\underline
{k}}(t,
\cdot) \bigr\|_{\mathscr{H}^{\otimes\underline{k}}}^2
\nonumber
\\[-8pt]
\label{norm-decomp}
\\[-8pt]
\nonumber
&&\quad= \bigl\| F^{(n_1)}_{\underline{k}}(t,\cdot) \bigr\|_{\mathscr{H}^{\otimes
\underline{k}}}^2
+ \bigl\|F^{(n_2)}_{\underline{k}}(t,\cdot) \bigr\| _{\mathscr{H}^{\otimes\underline{k}}}^2 -
2 \bigl\langle F^{(n_1)}_{\underline{k}}(t,\cdot), F^{(n_2)}_{\underline
{k}}(t,
\cdot ) \bigr\rangle_{\mathscr{H}^{\otimes\underline{k}}}.
\end{eqnarray}
By Definition
\eqref{F-def}, we have
\begin{eqnarray*}
&& \bigl\langle
F^{(n_1)}_{\underline{k}}(t,\cdot), F^{(n_2)}_{\underline
{k}}(t,
\cdot) \bigr\rangle_{\mathscr{H}^{\otimes\underline{k}}}
\\
&&\quad= a^2_{\underline{k}}\bigl\langle m(n_1) \bigr
\rangle^{\underline{k}-1}\bigl\langle m(n_2) \bigr\rangle^{\underline{k}-1}
\sum_{1 \leq i^{(1)} \leq\lfloor
m(n_1) t\rfloor} \sum_{1 \leq i^{(2)} \leq\lfloor m(n_2) t\rfloor
}
\bigl\langle g^{(n_1)}_{i^{(1)}},g^{(n_2)}_{i^{(2)}}
\bigr\rangle ^{\underline
{k}}_{\mathscr{H}^{\otimes\underline{k}}},
\end{eqnarray*}
where
$g^{(n)}_i$ is given by \eqref{hgdef}. Mimicking the proof
of Lemma~\ref{corr-struct}, we obtain %
\begin{eqnarray*}
&& \bigl\langle g^{(n_1)}_{i^{(1)}},g^{(n_2)}_{i^{(2)}}
\bigr\rangle_{\mathscr{H}^{\otimes\underline{k}}}
\\
&&\quad= \prod_{\nu=1}^d
\int G^{(1)}_{H_\nu} \biggl( \biggl[\frac
{i^{(1)}_\nu
-1}{m_\nu(n_1)},
\frac{i^{(1)}_\nu}{m_\nu(n_1)} \biggr), v \biggr)G^{(1)}_{H_\nu}
\biggl( \biggl[ \frac{i^{(2)}_\nu-1}{m_\nu
(n_2)},\frac
{i^{(2)}_\nu}{m_\nu(n_2)} \biggr), v \biggr)\, \mathrm{d}v
\\
&&\quad = \prod_{\nu=1}^d \mathbf{E} \biggl[
B_{H_\nu} \biggl( \biggl[\frac
{i^{(1)}_\nu
-1}{m_\nu(n_1)},\frac{i^{(1)}_\nu}{m_\nu(n_1)}
\biggr) \biggr) B_{H_\nu
} \biggl( \biggl[\frac{i^{(2)}_\nu-1}{m_\nu(n_2)},
\frac{i^{(2)}_\nu
}{m_\nu
(n_2)} \biggr) \biggr)\biggr]
\\
&&\quad = \prod_{\nu=1}^d H_\nu(2H_\nu-1)
\int_{({i^{(1)}_\nu -1})/{m_\nu
(n_1)}}^{{i^{(1)}_\nu}/{m_\nu(n_1)}}\int_{({i^{(2)}_\nu-1})/{m_\nu(n_2)}}^{{i^{(2)}_\nu}/{m_\nu(n_2)}}
|v_1-v_2|^{-2(1-H_\nu)} \,\mathrm{d} v_1
\,\mathrm{d}v_2,
\end{eqnarray*}
where the final equality follows (see, e.g., \cite{NVV1999}, page~574) since
$H_\nu> 1- \frac{1}{2\underline{k}}>\frac{1}{2}$ for any $\nu \in\{1,
\ldots,d\}$. Adapting the argument used in \cite{NNT}, pages~1064--1065, we deduce that
\begin{eqnarray}
&& \lim_{n_1,n_2 \rightarrow\infty} \bigl
\langle F^{(n_1)}_{\underline
{k}}(t,\cdot), F^{(n_2)}_{\underline{k}}(t,
\cdot) \bigr\rangle _{\mathscr
{H}^{\otimes\underline{k}}}
\nonumber
\\
&&\quad = a^2_{\underline{k}} \prod
_{\nu= 1}^d H_\nu^{\underline
{k}}(2H_\nu
-1)^{\underline{k}} \int_0^t\!\int
_0^t |v_1 - v_2|^{-2\underline
{k}(1-H_\nu)}
\,\mathrm{d}v_1 \,\mathrm{d}v_2
\nonumber
\\[-8pt]
\label{kernel-inner-product}
\\[-8pt]
\nonumber
&&\quad = a^2_{\underline{k}} \prod_{\nu= 1}^d
t^{2\widetilde{H}_\nu
}_\nu H_\nu^{\underline{k}}(2H_\nu-1)^{\underline{k}}
\int_0^1\!\int_0^1
|v_1 - v_2|^{-2\underline{k}(1-H_\nu)}\,\mathrm{d}v_1
\,\mathrm{d}v_2
\\
&&\quad = a^2_{\underline{k}} \prod_{\nu= 1}^d
t^{2\widetilde{H}_\nu
}_\nu \kappa(H_\nu, \underline{k}) =
a^2_{\underline{k}}\bigl\langle b^{(\underline{k})}\bigr\rangle
R^{(d)}_{\widetilde{H}}(t,t).
\nonumber
\end{eqnarray}
Thus, by \eqref{cov-conv} and \eqref{norm-decomp},
\begin{eqnarray*}
&& \lim_{n_1,n_2 \rightarrow\infty} \bigl\|F^{(n_1)}_{\underline
{k}}(t,
\cdot)-F^{(n_2)}_{\underline{k}}(t,\cdot) \bigr\|_{\mathscr
{H}^{\otimes\underline{k}}}^2
= 0,
\end{eqnarray*}
whence $ (F^{(n)}_{\underline{k}}(t,\cdot) )_{n \in\mathbb
{N}}$ is a
Cauchy sequence.

To characterize the limit of $ (F^{(n)}_{\underline{k}}(t,\cdot
)
)_{n \in\mathbb{N}}$, let us consider for any $s^{(1)},\ldots
,s^{(\underline
{k})}\in\mathbb{R}^d$,
\begin{eqnarray*}
&& F^{(n)}_{\underline{k}} \bigl(t,s^{(1)},
\ldots,s^{(\underline
{k})} \bigr) \\
&&\quad = a_{\underline{k}} \bigl\langle m(n) \bigr
\rangle^{\underline{k}-1} \sum_{1
\leq
i \leq\lfloor m(n) t \rfloor} \prod
_{\kappa=1}^{\underline{k}} G^{(d)}_H \biggl( \biggl[
\frac{i-1}{m(n)},\frac{i}{m(n)} \biggr),s^{(\kappa
)} \biggr)
\\
&&\quad = a_{\underline{k}} \bigl\langle m(n) \bigr\rangle^{\underline{k}-1} \sum
_{1
\leq i \leq\lfloor m(n) t \rfloor} \prod_{\kappa=1}^{\underline{k}}
\prod_{\nu=1}^d G^{(1)}_{H_\nu}
\biggl( \biggl[\frac{i_\nu
-1}{m_\nu
(n)},\frac{i_\nu}{m_\nu(n)} \biggr),s_\nu^{(\kappa)}
\biggr)
\\
&&\quad = a_{\underline{k}} \prod_{\nu=1}^d
\frac{1}{m_\nu(n)}\sum_{j=1}^{
\lfloor m_\nu(n) t_\nu\rfloor} \prod
_{\kappa=1}^{\underline{k}} m_\nu
(n)G^{(1)}_{H_\nu} \biggl( \biggl[\frac{j-1}{m_\nu(n)},
\frac
{j}{m_\nu
(n)} \biggr),s_\nu^{(\kappa)} \biggr),
\end{eqnarray*}
where the second equality is a consequence of Remark~\ref{tensorcase}. Since
\begin{eqnarray*}
&& G^{(1)}_{H_\nu} \biggl( \biggl[\frac{j-1}{m_\nu(n)},\frac{j}{m_\nu
(n)}
\biggr),s_\nu^{(\kappa)} \biggr)\\
&&\quad = \frac{1}{\chi
(H_\nu
)}
\biggl( \biggl(\frac{j}{m_\nu(n)} - s_\nu^{(\kappa)}
\biggr)_+^{H_\nu-
{1}/{2}} - \biggl(\frac{j-1}{m_\nu(n)} - s_\nu^{(\kappa)}
\biggr)_+^{H_\nu- {1}/{2}} \biggr),
\end{eqnarray*}
it follows from Lemma~\ref{integral-lemma}, below, that
\begin{equation}
\label{kernel-conv-constant}
F^{(n)}_{\underline{k}}(t,\cdot) \mathop{\longrightarrow}_{n
\rightarrow\infty} C'''
(a_{\underline{k}},H,\underline{k} ) \widehat {G}^{(\underline
{k})}_{\widetilde{H}}
(t,\cdot) \qquad \mbox{a.e. on $\mathbb{R} ^{\underline{k}d}$}
\end{equation}
for some constant $C''' (a_{\underline{k}},H,\underline{k} )>0$.
By the Cauchy property of $ (F^{(n)}_{\underline{k}}(t,\cdot)
)_{n \in\mathbb{N}}$, the convergence \eqref{kernel-conv-constant}
holds also
in $\mathscr{H}^{\otimes\underline{k}}$.
Clarke De la Cerda and Tudor \cite{CT2012}, pages~4--6,  have shown that
$\mathbf{E} [\widehat{Z}(t)^2 ] = \underline{k}! \|
\widehat
{G}^{(\underline{k})}_{\widetilde{H}} (t,\cdot) \|^2_{\mathscr
{H}^{\otimes\underline{k}}} = R^{(d)}_{\widetilde{H}}(t,t)$. In view
of \eqref{kernel-inner-product}, we find that
\begin{eqnarray*}
&& C''' (a_{\underline{k}},H,\underline{k}
)^2 = \underline{k}! a^2_{\underline{k}} \bigl\langle
b^{(\underline{k})} \bigr\rangle= \Lambda_{H,f},
\end{eqnarray*}
whence \eqref{kernel-convergence} follows.
\end{pf}

The following technical lemma was essential in the proof of Proposition~\ref{PNCLT}.

\begin{lem}\label{integral-lemma} Suppose that $k \geq2$, $\check{H}
\in (\frac{1}{2},1 )$, and $v>0$. Then
\begin{eqnarray}
&& \frac{1}{n} \sum_{j=1}^{\lfloor nv \rfloor}
\prod_{\kappa=1}^{k} n \biggl( \biggl(
\frac{j}{n} - s_\kappa \biggr)_+^{\check{H}-{1}/{2}}- \biggl(
\frac{j-1}{n} - s_\kappa \biggr)_+^{\check{H}-{1}/{2}} \biggr)
\nonumber
\\[-8pt]
\label{Hermite-limit}
\\[-8pt]
\nonumber
&&\quad\mathop{\longrightarrow}_{n \rightarrow\infty} \biggl(\check{H}-\frac
{1}{2}
\biggr)^k \int_0^v \prod
_{\kappa=1}^k (u-s_\kappa)_+^{\check{H}-{3}/{2}}
\,\mathrm{d}u
\end{eqnarray}
for almost any $s=(s_1,\ldots,s_k) \in\mathbb{R}^k$.
\end{lem}

\begin{pf}
We may assume that $\overline{s} :=\max(s_1,\ldots,s_k) < v$, as
otherwise \eqref{Hermite-limit} is trivially true. In fact,
\begin{eqnarray*}
&& \int_0^v \prod
_{\kappa=1}^k (y-s_\kappa)_+^{\check{H}-{3}/{2}}
\,\mathrm{d} y = \int_{\overline{s}}^v \prod
_{\kappa=1}^k (y-s_\kappa)^{\check
{H}-{3}/{2}}
\,\mathrm{d}y.
\end{eqnarray*}
We split the sum on the left-hand side of \eqref{Hermite-limit} for any
$n \in\mathbb{N}$, such that $\lfloor nv \rfloor> \lfloor n\overline
{s}\rfloor+ 3$, as
\begin{eqnarray*}
&& \frac{1}{n} \sum_{j=1}^{\lfloor nv \rfloor} \prod
_{\kappa=1}^{k} n \biggl( \biggl(
\frac{j}{n} - s_\kappa \biggr)_+^{\check{H}-{1}/{2}}- \biggl(
\frac{j-1}{n} - s_\kappa \biggr)_+^{\check{H}-{1}/{2}} \biggr)
\\
&&\quad = \frac{1}{n} \sum
_{j=\lfloor n\overline{s} \rfloor+1}^{\lfloor
n\overline{s} \rfloor+ 2} \prod_{\kappa=1}^{k}
n \biggl( \biggl( \frac
{j}{n} - s_\kappa \biggr)^{\check{H}-{1}/{2}}-
\biggl( \frac
{j-1}{n} - s_\kappa \biggr)_+^{\check{H}-{1}/{2}} \biggr)
\\
&&\qquad  {}+ \frac{1}{n} \sum_{j=\lfloor n\overline{s} \rfloor
+3}^{\lfloor nv \rfloor}
\prod_{\kappa=1}^{k} n \biggl( \biggl(
\frac
{j}{n} - s_\kappa \biggr)^{\check{H}-{1}/{2}}- \biggl(
\frac
{j-1}{n} - s_\kappa \biggr)^{\check{H}-{1}/{2}} \biggr)
\\
&&\quad =:S^{(1)}_n + S^{(2)}_n.
\end{eqnarray*}
Using the mean value theorem, we obtain for any $y \in\mathbb{R}$ and
$n$,
$j\in\mathbb{N}$, such that $\frac{j-1}{n}>y$, the bounds
\begin{eqnarray}\label{mv-upper}
n \biggl( \biggl( \frac{j}{n} - y \biggr)^{\check{H}-{1}/{2}}- \biggl(
\frac{j-1}{n} - y \biggr)^{\check{H}-{1}/{2}} \biggr) &\leq &  \biggl(\check {H}-
\frac{1}{2} \biggr) \biggl( \frac{j-1}{n}- y \biggr)^{\check
{H}-{3}/{2}},
\\
\label{mv-lower}
n \biggl( \biggl( \frac{j}{n} - y \biggr)^{\check{H}-{1}/{2}}- \biggl(
\frac{j-1}{n} - y \biggr)^{\check{H}-{1}/{2}} \biggr) & \geq &  \biggl(\check{H}-
\frac{1}{2} \biggr) \biggl( \frac{j}{n}- y \biggr)^{\check
{H}-{3}/{2}}.
\end{eqnarray}
Since we are aiming to prove \eqref{Hermite-limit} for \emph{almost}
any $s \in\mathbb{R}^k$, we may assume (by symmetry) that $\overline
{s} =
s_1>s_\kappa$ for any $\kappa\in\{2,\ldots,k\}$. Then we have for $j
\in\{1,2\}$,
\begin{eqnarray*}
&& \limsup_{n \rightarrow\infty}\prod_{\kappa=2}^{k}
n \biggl( \biggl( \frac
{\lfloor n\overline{s} \rfloor+ j}{n} - s_\kappa \biggr)^{\check
{H}-{1}/{2}}-
\biggl( \frac{\lfloor n\overline{s}\rfloor+j-1}{n} - s_\kappa \biggr)_+^{\check{H}-{1}/{2}} \biggr)
<\infty
\end{eqnarray*}
by \eqref{mv-upper}, and
\begin{eqnarray*}
&& 0 \leq \biggl( \frac{\lfloor n\overline{s} \rfloor+ j}{n} - s_1 \biggr)^{\check{H}-{1}/{2}}-
\biggl( \frac{\lfloor n\overline
{s}\rfloor
+j-1}{n} - s_1 \biggr)_+^{\check{H}-{1}/{2}}
\leq \biggl(\frac
{\lfloor n\overline{s} \rfloor+ 2}{n} - s_1 \biggr)^{\check{H}-{1}/{2}}
\mathop{\longrightarrow}_{n \rightarrow\infty} 0.
\end{eqnarray*}
Hence, we find that $S^{(1)}_n \rightarrow0$ as $n \rightarrow\infty$.

Finally, invoking \eqref{mv-upper}, we obtain
\begin{eqnarray*}
S^{(2)}_n & \leq & \biggl(
\check{H}-\frac{1}{2} \biggr)^k \frac{1}{n} \sum
_{j=\lfloor n \overline{s}\rfloor+ 3}^{\lfloor nv\rfloor} \prod_{\kappa= 1}^k
\biggl( \frac{j-1}{n} -s_\kappa \biggr)^{\check
{H}-{3}/{2}}
\\
& =& \biggl( \check{H}-\frac{1}{2} \biggr)^k \int
_{({\lfloor n
\overline{s}\rfloor+ 2})/{n}}^{{\lfloor nv \rfloor}/{n}} \prod_{\kappa= 1}^k
\biggl(\frac{ \lfloor n y \rfloor+ 1}{n} - \frac{1}{n} - s_\kappa
\biggr)^{\check{H}-{3}/{2}} \,\mathrm{d}y
\\
& \leq &  \biggl( \check{H}-\frac{1}{2} \biggr)^k \int
_{({\lfloor n \overline{s}\rfloor+ 1})/{n}}^{({\lfloor nv \rfloor-1})/{n}} \prod_{\kappa= 1}^k
(y-s_\kappa)^{\check{H}- {3}/{2}} \,\mathrm{d}y \mathop{\longrightarrow}_{n \rightarrow\infty} \biggl( \widehat{H}-\frac{1}{2} \biggr)^k \int
_{
\overline{s}}^{v} \prod_{\kappa= 1}^k
(y-s_\kappa)^{\check{H}-
{3}/{2}} \,\mathrm{d}y
\end{eqnarray*}
and similarly by \eqref{mv-lower},
\begin{eqnarray*}
S^{(2)}_n & \geq &  \biggl(
\check{H}-\frac{1}{2} \biggr)^k \frac{1}{n} \sum
_{j=\lfloor n \overline{s}\rfloor+ 3}^{\lfloor nv\rfloor} \prod_{\kappa= 1}^k
\biggl( \frac{j}{n} -s_\kappa \biggr)^{\check
{H}-{3}/{2}}
\\
& =& \biggl( \check{H}-\frac{1}{2} \biggr)^k \int
_{({\lfloor n
\overline
{s}\rfloor+ 2})/{n}}^{{\lfloor nv \rfloor}/{n}} \prod_{\kappa= 1}^k
\biggl(\frac{ \lfloor n y \rfloor+ 1}{n} - s_\kappa \biggr)^{\check
{H}-{3}/{2}}
\,\mathrm{d}y
\\
&\geq & \biggl( \check{H}-\frac{1}{2} \biggr)^k \int
_{({\lfloor n \overline{s}\rfloor+ 3})/{n}}^{({\lfloor nv \rfloor+1})/{n}} \prod_{\kappa= 1}^k
(y-s_\kappa)^{\check{H}- {3}/{2}} \,\mathrm{d}y \mathop{\longrightarrow}_{n \rightarrow
\infty} \biggl( \check{H}-\frac{1}{2} \biggr)^k \int
_{
\overline{s}}^{v} \prod_{\kappa= 1}^k
(y-s_\kappa)^{\check{H}-
{3}/{2}} \,\mathrm{d}y.
\end{eqnarray*}
(The convergence of the bounding integrals above, as $n \rightarrow
\infty$, is ensured by Lebesgue's dominated convergence theorem.) Thus,
the convergence \eqref{Hermite-limit} follows from the sandwich
lemma.\vfill
\end{pf}

\section{Functional convergence} \label{fun-convergence}

To show that Theorems \ref{gvar-clt} and \ref{gvar-nclt} indeed hold in
the functional sense, we need to establish tightness of the relevant
families of processes in the space $D([0,1]^d)$. To this end, we use
the tightness criterion due to Bickel and Wichura \cite{BW1971}, Theorem~3. To apply this criterion, we need to bound the fourth
moments of the increments of $\overline{U}^{(n)}_f$ uniformly over $n
\in\mathbb{N}$.

\subsection{Moment bound and diagrams}

As a preparation for the proof of tightness, we establish a moment
bound for nonlinear functionals of stationary Gaussian processes
indexed by $\mathbb{N}^d$.
The bound is a multi-parameter extension of Proposition~4.2 of \cite
{Taqqu1977}, albeit under more restrictive assumptions.

\begin{lem}[(Moment bound)]\label{moment-bound}
Let $f$ be as in Section~\ref{secmainresults} and $\{Y_{i} \dvt i \in\mathbb
{N}^d\}$
a Gaussian process such that $\mathbf{E}[Y_i]=0$ and $\mathbf
{E}[Y_i^2]=1$ for any $i
\in\mathbb{N}^d$. Moreover, suppose that there exists a function
$\rho\dvt \mathbb{Z}^d
\rightarrow[-1,1]$ such that $\mathbf{E}[Y_{i^{(1)}} Y_{i^{(2)}}]=
\rho
(i^{(1)}-i^{(2)})$ for any $i^{(1)}$, $i^{(2)} \in\mathbb{N}^d$. If
$p \in\{
2,3,\ldots\}$ and the Hermite coefficients $a_{\underline
{k}},a_{\underline{k}+1},\ldots$ of the function $f$ satisfy
\begin{eqnarray*}
&& C''''(f,p) :=\sum
_{k=\underline{k}}^{\infty} (p-1)^{k/2} \sqrt {k!}
|a_k|<\infty,
\end{eqnarray*}
then for any $l \in\mathbb{N}^d$,
\begin{eqnarray*}
&& \biggl|\mathbf{E} \biggl[ \biggl(\langle l\rangle^{-1/2} \sum
_{1 \leq
i \leq l} f(Y_{i}) \biggr)^p \biggr] \biggr|
\leq \biggl(2^d C''''(f,p)^2
\sum_{|i|<l} \bigl|\rho(i)\bigr|^{\underline{k}}
\biggr)^{p/2}.
\end{eqnarray*}
\end{lem}

The proof of Proposition~4.2 of \cite{Taqqu1977} is based on a graph
theoretic argument that involves \emph{multi-graphs}. We prove Lemma~\ref{moment-bound} using slightly different (but essentially analogous)
formalism based on \emph{diagrams}, defined below. Breuer and Major
\cite{BM} used diagrams to prove their central limit theorem for
nonlinear functionals of Gaussian random fields via the \emph{method of
moments}. In fact in the proof of Lemma~\ref{moment-bound}, we adapt
some of the arguments used in \cite{BM}.

\begin{defn}
Let $p \in\{2,3,\ldots\}$ and $(k_1,\ldots,k_p) \in\mathbb{N}^p$
be such that
$k_1+\cdots+k_p$ is an even number.
A \emph{diagram} of order $(k_1,\ldots,k_p)$ is a graph $G=(V_G,E_G)$
with the following three properties:
\begin{enumerate}[(3)]
\item[(1)] We have ${  V_G=\bigcup_{j=1}^p \{(j,1),\ldots
,(j,k_{j})\}}$.
\item[(2)] The degree of any vertex $v \in V_G$ is one.
\item[(3)] Any edge $ ((j,k),(j',k') ) \in E_G$ has the property that
$j \neq j'$.
\end{enumerate}
We denote the class of diagrams of order $(k_1,\ldots,k_p)$ by
$\mathscr
{G}(k_1,\ldots,k_p)$.
For the sake of completeness we set $\mathscr{G}(k_1,\ldots,k_p)
:=\varnothing$ when $k_1+\cdots+k_p$ is an odd number (no
diagrams can then exist by the \emph{handshaking lemma} of graph
theory). Let us also define two functions $\lambda_1$ and $\lambda_2$
of an edge $e= ((j,k),(j',k') )\in E_G$, where $j<j'$, by setting
$\lambda_1(e) :=j$ and $\lambda_2(e) :=j'$.
\end{defn}

Diagrams are connected to Hermite polynomials and Gaussian random
variables via the so-called \emph{diagram formula}, which is originally
due to Taqqu \cite{Taqqu1977}, Lemma~3.2. Below, we state a version of
the formula that appears in \cite{BM}, page~431.

\begin{lem}[(Diagram formula)]\label{diagram}
Let $p \in\{2,3,\ldots\}$ and let $Y_1,\ldots,Y_p$ be jointly Gaussian
random variables with $\mathbf{E}[Y_i]=0$ and $\mathbf{E}[Y^2_i]=1$
for any $i\in\{
1,\ldots,p\}$.
For any $(k_1,\ldots,k_p) \in\mathbb{N}^p$, we have\vspace*{-8pt}
\begin{eqnarray*}
&& \mathbf{E} \Biggl[\prod_{j=1}^p
P_{k_j}(Y_j) \Biggr] = \sum_{G \in
\mathscr
{G}(k_1,\ldots,k_p)}
\prod_{e \in E_G} \mathbf{E} [Y_{\lambda_1(e)}Y_{\lambda_2(e)}
],
\end{eqnarray*}
where a sum over an empty index set is interpreted as\vspace*{-3pt} zero.
\end{lem}

\begin{rem}
The diagram formula can be used to estimate the cardinalities of
classes of diagrams.
As pointed out by Bardet and Surgailis \cite{BS2013}, page~461,
using Lemma~\ref{diagram} and Lemma~3.1 of \cite{Taqqu1977} in the
special case $Y :=Y_1 = \cdots= Y_p$, we\vspace*{-3pt} obtain
\begin{equation}
\label{diagrambound}
\bigl|\mathscr{G}(k_1,\ldots,k_p)\bigr| =
\mathbf{E} \Biggl[\prod_{j=1}^p
P_{k_j}(Y) \Biggr]\leq(p-1)^{(k_1+\cdots+k_p)/2} \sqrt{k_1 !
\cdots k_p !}.
\end{equation}
\end{rem}

\begin{pf*}{Proof of Lemma~\ref{moment-bound}}
Fix $l \in\mathbb
{N}^d$. Let
us define for any $K\geq\underline{k}$, a polynomial\vspace*{-7pt} function
\begin{eqnarray*}
&& f_K(x) = \sum_{k=\underline{k}}^K
a_k P_k(x),   \qquad x \in\mathbb{R}.
\end{eqnarray*}
By Fatou's lemma, Lemma~\ref{diagram} and inequality \eqref
{diagrambound}, it follows\vspace*{-4pt} that
\begin{eqnarray*}
\mathbf{E} \bigl[ \bigl|f(Y_{i})-f_K(Y_{i})
\bigr|^p \bigr] & \leq & \sum_{k_1,\ldots,k_p
= K+1}^{\infty}
|a_{k_1}\cdots a_{k_p}| \bigl|\mathscr{G}(k_1,\ldots
,k_p)\bigr|
\\[-4pt]
& \leq &  \Biggl( \sum_{k = K+1}^\infty(p-1)^{k/2}
\sqrt{k!}|a_k| \Biggr)^p \mathop{\longrightarrow}_{K \rightarrow
\infty} 0
\end{eqnarray*}
for any $i \in\mathbb{N}^d$.
Thus, if $\varepsilon>0$, then there exists $K(l) \in\mathbb{N}$
such\vspace*{-4pt} that
\begin{equation}
\label{truncation}
\biggl|\mathbf{E} \biggl[ \biggl(\langle l\rangle^{-1/2} \sum
_{1 \leq
i \leq l} f(Y_{i}) \biggr)^p
\biggr]-\mathbf{E} \biggl[ \biggl(\langle l\rangle^{-1/2} \sum
_{1
\leq i \leq l} f_{K(l)}(Y_{i})
\biggr)^p \biggr] \biggr| \leq\varepsilon,
\end{equation}
by Minkowski's inequality and the fact that $X_n \displaystyle{\mathop{\longrightarrow}^{L^p(\Omega)}} X$ implies $\mathbf{E}[X^p_n] \rightarrow\mathbf
{E}[X^p]$ when $p
\in
\{2,3,\ldots\}$.
Lemma~\ref{diagram} yields now the\vspace*{-4pt} expansion
\begin{eqnarray*}
&& \mathbf{E} \biggl[ \biggl(\langle l\rangle^{-1/2} \sum
_{1 \leq i \leq l} f_{K(l)}(Y_{i})
\biggr)^p \biggr]
\\[-4pt]
&&\quad = \langle l\rangle^{-p/2}\mathop{\sum_{1 \leq i^{(j)}\leq l}}_{j
\in
\{1,\ldots,p\}}\sum_{k_1,\ldots,k_p = \underline{k}}^{K(l)}
a_{k_1}\cdots a_{k_p} \sum_{G \in\mathscr{G}(k_1,\ldots,k_p)}
\prod_{e
\in E_G} \mathbf{E} [Y_{i^{(\lambda_1(e))}}Y_{i^{(\lambda
_2(e))}}
]
\\[-4pt]
&&\quad = \sum_{k_1,\ldots,k_p = \underline{k}}^{K(l)} a_{k_1}
\cdots a_{k_p} \sum_{G \in\mathscr{G}(k_1,\ldots,k_p)}
I_G(l),
\end{eqnarray*}
where
\begin{equation}
\label{I-def}
I_G(l) :=\langle l\rangle^{-p/2}\mathop{\sum_{1 \leq
i^{(j)}\leq l}}_{j \in\{1,\ldots,p\}} \prod_{e \in E_G} \rho
\bigl(i^{(\lambda_1(e))}-i^{(\lambda_2(e))} \bigr), \qquad  G \in\mathscr
{G}(k_1,\ldots,k_p).
\end{equation}
By Lemma~\ref{diagramcor} below and inequality \eqref{diagrambound}, we
obtain the bound
\begin{eqnarray*}
&& \Biggl|\sum_{k_1,\ldots,k_p = \underline{k}}^{K(l)} a_{k_1}\cdots
a_{k_p} \sum_{G \in\mathscr{G}(k_1,\ldots,k_p)} I_G(l) \Biggr|
\\
&&\quad\leq \Biggl(\sum_{k = \underline{k}}^{K(l)}
(p-1)^{k/2}\sqrt {k!}|a_k| \Biggr)^p
\biggl(2^d\sum_{|i|<l} \bigl|
\rho(i)\bigr|^{\underline{k}} \biggr)^{p/2} \leq \biggl(2^d
C''''(f,p)^2
\sum_{|i|<l} \bigl|\rho(i)\bigr|^{\underline{k}}
\biggr)^{p/2}.
\end{eqnarray*}
In view of \eqref{truncation},
\begin{eqnarray*}
&&\biggl|\mathbf{E} \biggl[ \biggl(\langle l\rangle^{-1/2} \sum
_{1 \leq
i \leq l} f(Y_{i}) \biggr)^p \biggr] \biggr|
\leq \biggl(2^dC''''(f,p)^2
\sum_{|i|<l} \bigl|\rho(i)\bigr|^{\underline{k}}
\biggr)^{p/2} + \varepsilon,
\end{eqnarray*}
and letting $\varepsilon\rightarrow0$ concludes the proof.
\end{pf*}

The key ingredient in the proof of Lemma~\ref{moment-bound} was the
following uniform bound for the absolute value of the quantity
$I_G(l)$. We will derive this bound by adapting the asymptotic analysis
of the moments of a nonlinear functional of a Gaussian random field,
carried out in \cite{BM}, pages~435--436.

\begin{lem}\label{diagramcor}
For any $k_1,\ldots,k_p \geq\underline
{k}$, $G \in\mathscr{G}(k_1,\ldots,k_p)$, and $l \in\mathbb{N}^d$,
\begin{eqnarray*}
&& \bigl|I_G(l)\bigr| \leq \biggl(2^d\sum
_{|i|<l} \bigl|\rho(i)\bigr|^{\underline{k}} \biggr)^{p/2},
\end{eqnarray*}
where $I_G(l)$ is defined by \eqref{I-def}.
\end{lem}

\begin{pf}
As pointed out by Breuer and Major \cite{BM}, page~435, the quantity
$I_G(l)$ is invariant under permutation of the levels of the diagram
$G$. More precisely, if $\sigma$ is a permutation of the set $\{
1,\ldots
,p\}$, then we define a new diagram $\widetilde{G} \in\mathscr
{G}(k_{\sigma
(1)},\ldots,k_{\sigma(p)})$ such that $ ((j,k),(j',k') ) \in
E_{\widetilde{G}}$ if and only if $ ((\sigma^{-1}(j),k),(\sigma
^{-1}(j'),k') ) \in E_{G}$. For such a diagram $\widetilde{G}$ it holds
that $I_G(l) = I_{\widetilde{G}}(l)$. Relying on this invariance property
we assume, without loss of generality, that
\begin{equation}
\label{incrcard} k_1 \leq k_2 \leq\cdots\leq
k_{p-1} \leq k_p.
\end{equation}

Let us introduce the notation $k_G(j) :=|\{ e \in E_G \dvt \lambda
_1(e) = j \}| \in\{0,1,\ldots,k_{j}\}$ for any $j \in\{1,\ldots,p\}$.
Since $\lambda_1(e)<\lambda_2(e)$ for any $e \in E_G$, we have
\begin{eqnarray}
\bigl|I_G(l)\bigr| & \leq & \langle l
\rangle^{-p/2} \mathop{\sum_{1 \leq i^{(\kappa)}\leq l}}_{\kappa\in\{1,\ldots,p\}} \prod
_{j=1}^p \mathop{\prod_{e \in E_G}}_{\lambda_1(e) = j}
\bigl|\rho \bigl(i^{(j)}-i^{(\lambda_2(e))} \bigr) \bigr|
\nonumber
\\[-8pt]
\label{separation}
\\[-8pt]
\nonumber
& =&  \langle l\rangle^{-p/2}
\mathop{\sum_{1 \leq i^{(\kappa)}\leq l}}_{\kappa\in\{2,\ldots,p\}}
\prod_{j=2}^p \mathop{\prod_{e \in E_G}}_{\lambda_1(e) = j} \bigl|\rho \bigl(i^{(j)}-i^{(\lambda_2(e))} \bigr) \bigr|
\sum_{1 \leq i^{(1)}\leq l}\mathop{\prod_{e \in E_G}}_{\lambda_1(e)
= 1} \bigl|
\rho \bigl(i^{(1)}-i^{(\lambda_2(e))} \bigr) \bigr|.
\end{eqnarray}
Using Young's inequality (see \cite{BM}, page~435) and the trivial estimate
\begin{eqnarray*}
&& \sup_{1 \leq i \leq l}\sum_{1 \leq i^{(1)}\leq l} \bigl|\rho
\bigl(i^{(1)}-i \bigr) \bigr|^q \leq\sum
_{|i|<l} \bigl|\rho(i)\bigr|^q,  \qquad q \geq0,
\end{eqnarray*}
one can show that
\begin{eqnarray*}
&& \mathop{\sup_{1 \leq i^{(\kappa)}\leq l}}_{\kappa\in\{2,\ldots
,p\}}\sum_{1 \leq i^{(1)}\leq l}\mathop{\prod_{e \in E_G}}_{\lambda_1(e)
= 1} \bigl|\rho \bigl(i^{(1)}-i^{(\lambda_2(e))}
\bigr) \bigr| \leq \sum_{|i|<l} \bigl|\rho(i)\bigr|^{k_G(1)}.
\end{eqnarray*}
Applying this procedure, mutatis mutandis, to \eqref{separation}
repeatedly we arrive at
\begin{equation}
\label{young}
\bigl|I_G(l)\bigr| \leq\langle l\rangle^{-p/2} \prod
_{j=1}^p \sum
_{|i|< l} \bigl|\rho (i)\bigr|^{k_G(j)}.
\end{equation}

By H\"older's inequality, we have for any $j \in\{1,\ldots,p\}$,
\begin{eqnarray*}
\sum_{|i|<l} \bigl|
\rho(i)\bigr|^{k_G(j)} & \leq & \langle2 l \rangle^{1-
k_G(j)/k_j} \biggl(\sum
_{|i|<l} \bigl|\rho(i)\bigr|^{k_j}
\biggr)^{k_G(j)/k_j}
\\
& \leq & \langle2 l \rangle^{1- k_G(j)/k_j} \biggl(\sum
_{|j|<l} \bigl|\rho (i)\bigr|^{\underline{k}} \biggr)^{k_G(j)/k_j},
\end{eqnarray*}
where we use the proviso $k_j \geq\underline{k}$ to deduce the second
inequality. Returning to \eqref{young}, we have thus established that
\begin{equation}
\label{I-bound}
\bigl|I_G(l)\bigr| \leq \bigl(2^d
\bigr)^{p/2}\langle2 l \rangle^{p/2- \sum_{j=1}^p
k_G(j)/k_j} \biggl(\sum
_{|i|<l} \bigl|\rho(i)\bigr|^{\underline{k}} \biggr)^{\sum
_{j=1}^p k_G(j)/k_j}.
\end{equation}
Breuer and Major \cite{BM}, page 436, have shown that whenever \eqref
{incrcard} holds, we have
\begin{equation}
\label{levelbound}
\sum_{j=1}^p
\frac{k_G(j)}{k_j} - \frac{p}{2} \geq0
\end{equation}
(see also Remark~\ref{non-regular}, below). By \eqref{levelbound}, we
may use the rough estimate $\sum_{|i|<l} |\rho(i)|^{\underline{k}}
\leq
\langle2 l\rangle$ to deduce that
\begin{eqnarray}
\biggl(\sum_{|i|<l }
\bigl|\rho(i)\bigr|^{\underline{k}} \biggr)^{\sum_{j=1}^p
k_G(j)/k_j} & =& \biggl(\sum
_{|i|<l} \bigl|\rho(i)\bigr|^{\underline{k}} \biggr)^{\sum
_{j=1}^p k_G(j)/k_j -p/2}
\biggl(\sum_{|i|<l} \bigl|\rho(i)\bigr|^{\underline
{k}}
\biggr)^{p/2}
\nonumber
\\[-8pt]
\label{rho-sum-bound} %
\\[-8pt]
\nonumber
& \leq & \langle2 l \rangle^{\sum_{j=1}^p k_G(j)/k_j -p/2} \biggl(\sum
_{|i|<l} \bigl|\rho(i)\bigr|^{\underline{k}} \biggr)^{p/2}.
\end{eqnarray}
The assertion follows now by applying \eqref{rho-sum-bound} to \eqref
{I-bound}.
\end{pf}

\begin{rem}\label{non-regular}
Strictly speaking, the inequality \eqref{levelbound} is shown in \cite
{BM} as a part of a more extensive argument that uses the assumption
that the diagram $G$ is not \emph{regular} (see \cite{BM}, page 432, for
the definition of regularity). However, the assumption of
non-regularity of $G$ is completely immaterial concerning the validity
of \eqref{levelbound} and, in fact, not used in the proof in \cite{BM}, page~436.
\end{rem}

\subsection{Tightness} \label{ssectightness}

Furnished with the moment bound of Lemma~\ref{moment-bound}, we prove
the following lemma that enables us to complete the proofs of Theorems
\ref{gvar-clt} and \ref{gvar-nclt}.

\begin{lem}[(Tightness)]\label{tightness-lem}
Suppose that $H \in(0,1)^d$ and that Assumption~\ref
{Hermite-assumption} holds. Then, the family $ \{\overline
{U}^{(n)}_f \dvt n \in\mathbb{N} \}$ is tight in $D([0,1]^d)$.
\end{lem}

\begin{pf} The assertion follows from Theorem~3 of \cite{BW1971},
provided that
\begin{equation}
\label{tightness-condition}
\sup_{n \in\mathbb{N}^d} \mathop{\sup_{s,  t\in[0,1]^d}}_{s <
t}
\frac
{\mathbf{E}
 [\overline{U}^{(n)}_f([s,t))^4 ]}{\langle t-s\rangle^2} < \infty.
\end{equation}
But since for any $n\in\mathbb{N}$, the realization of $\overline{U}^{(n)}_f$
is constant on any set of the form
\begin{eqnarray*}
&& \biggl[\frac{i-1}{m(n)},\frac{i}{m(n)} \biggr),\qquad  1 \leq i \leq m(n),
\end{eqnarray*}
it suffices to show (see \cite{BW1971}, page~1665) that
\begin{equation}
\label{tightness-condition-2}
\sup_{n \in\mathbb{N}} \mathop{\sup_{s, t \in\mathscr{E}_n}}_{s < t}
\frac
{\mathbf{E} [\overline{U}^{(n)}_f([s,t))^4 ]}{\langle
t-s\rangle
^2}<\infty,
\end{equation}
where $\mathscr{E}_n :=\{ i/m(n) \dvt 0 \leq i \leq m(n) \}$,
instead of \eqref{tightness-condition}.

Using Lemmas \ref{corr-struct} and \ref{moment-bound}, we arrive at
\begin{eqnarray*}
\sup_{n \in\mathbb{N}}\mathop{\sup_{s,t \in\mathscr{E}_n}}_{s
< t}
\frac
{\mathbf{E} [\overline{U}^{(n)}_f([s,t))^4 ]}{\langle
t-s\rangle^2} & =& \sup_{n \in\mathbb{N}} \biggl\langle
\frac{ m(n) }{ c(n)} \biggr\rangle ^2\sup_{1 \leq l \leq m(n)}
\mathbf{E} \biggl[ \biggl(\langle l\rangle^{-1/2} \sum
_{1\leq i \leq l} f \bigl( X^{(n)}_{i} \bigr)
\biggr)^4 \biggr]
\\
& \leq & \sup_{n \in\mathbb{N}} \Biggl(2^d C''''(f,4)
\prod_{\nu
=1}^d \psi _\nu (n)
\Biggr)^2,
\end{eqnarray*}
where
\begin{eqnarray*}
\psi_\nu(n) :=
\cases{ \displaystyle{
\frac{1}{m_\nu(n)^{1-2\underline{k}(1-H_\nu)}}\sum_{|j|<m_\nu(n)} \bigl|r_{H_\nu}(j)\bigr|^{\underline{k}}}
\leq C'(H_\nu ,\underline {k}), &\quad $\displaystyle H_\nu\in
\biggl(1-\frac{1}{2\underline{k}},1 \biggr)$,\vspace*{3pt}
\cr
\displaystyle{ \frac{1}{\log (m_\nu(n) )}\sum_{|j|<m_\nu(n)}
\bigl|r_{H_\nu}(j)\bigr|^{\underline{k}}} \leq C'(H_\nu,
\underline{k}), & \quad $\displaystyle H_\nu= 1-\frac{1}{2\underline{k}}$,\vspace*{3pt}
\cr
\displaystyle { \sum_{|j|<m_\nu(n)} \bigl|r_{H_\nu}(j)\bigr|^{\underline{k}}
\leq \sum_{j\in\mathbb{Z}} \bigl|r_{H_\nu}(j)\bigr|^{\underline{k}}}<
\infty, &\quad $\displaystyle H_\nu\in \biggl(0,1-\frac{1}{2\underline{k}} \biggr)$.}
\end{eqnarray*}
(The first two inequalities above follow from the estimate \eqref
{divergent-correlations}.)
We have, thus, verified the tightness condition \eqref{tightness-condition-2}.
\end{pf}

\begin{pf*}{Proof of Theorem~\ref{gvar-clt}}
Recall that, for a family of pairs of random elements, tightness of
marginals implies joint tightness. Thus, it follows from Lemma~\ref
{tightness-lem} that the family $ \{ (Z,  \overline
{U}^{(n)}_f ) \dvt n \in\mathbb{N} \}$ is tight in
$D([0,1]^d)^2$. The
assertion follows then from Proposition~\ref{fdd-conv-prop} and Theorem~2 of \cite{BW1971}.
\end{pf*}

\begin{pf*}{Proof of Theorem~\protect\ref{gvar-nclt}}
Analogously to the proof of Theorem~\ref{gvar-clt}, above, we deduce
from Lemma~\ref{tightness-lem} that $ \{ (\Lambda
_{H,f}^{{1}/{2}}\widehat{Z} , \overline{U}^{(n)}_f ) \dvt n \in\mathbb
{N} \}$ is
tight in $D([0,1]^d)^2$. Moreover, Proposition~\ref{PNCLT} implies that
\begin{eqnarray*}
&& \overline{U}^{(n)}_f(t) \mathop{\longrightarrow}_{n\rightarrow\infty}^{
\mathbf{P} }\Lambda _{H,f}^{{1}/{2}}\widehat{Z}(t),\qquad   t
\in[0,1]^d,
\end{eqnarray*}
which, in turn, implies the corresponding convergence of
finite-dimensional laws. Thus, by Theorem~2 of \cite{BW1971}, we have
\begin{equation}
\label{diagonal-conv}
\bigl(\Lambda_{H,f}^{{1}/{2}}\widehat{Z} ,
\overline {U}^{(n)}_f \bigr)\mathop{\longrightarrow}_{n\rightarrow\infty}^{\mathscr{L}} \bigl(\Lambda _{H,f}^{{1}/{2}}\widehat{Z} ,
\Lambda_{H,f}^{{1}/{2}}\widehat {Z} \bigr) \qquad \mbox{in }D
\bigl([0,1]^d\bigr)^2.
\end{equation}
Since the limit in \eqref{diagonal-conv} belongs to $C([0,1]^d)^2$ and
since substraction is a continuous operation on $C([0,1]^d)^2$ (with
respect to the Skorohod topology),
the continuous mapping theorem implies that
\begin{equation}
\label{difference}
\overline{U}^{(n)}_f -
\Lambda_{H,f}^{{1}/{2}}\widehat {Z}\mathop{\longrightarrow}_{n\rightarrow
\infty}^{\mathscr{L}} 0 \qquad \mbox{in }D\bigl([0,1]^d\bigr).
\end{equation}
It remains to note that the convergence \eqref{difference} holds also
in probability as the limit is deterministic.
\end{pf*}

\section{Application to power variations}\label{secpower-var}

\subsection{Convergence of power variations and their fluctuations}

As an application of Theorems \ref{gvar-clt} and \ref{gvar-nclt}, we
study the asymptotic behaviour of \emph{signed power variations} of the
fBs $Z$. Let $p \in\mathbb{N}$ be fixed throughout this section.
We consider a family $\{V^{(n)}_p \dvt n \in\mathbb{N}\}$ of $d$-parameter
processes, given by
\begin{eqnarray*}
&& V^{(n)}_{p}(t) :=\bigl\langle m(n)^{pH-1} \bigr
\rangle\sum_{1 \leq i \leq
\lfloor m(n)t \rfloor} Z \biggl( \biggl[\frac{i-1}{m(n)},
\frac
{i}{m(n)} \biggr) \biggr)^{p},\qquad   t \in[0,1]^d,
  n \in\mathbb{N}.
\end{eqnarray*}
The realizations of $V^{(n)}_{p}$ belong to the space $D([0,1]^d)$, as
was the case with generalized variations.
To describe the asymptotic behaviour of $V^{(n)}_{p}$, we introduce
\begin{eqnarray*}
v_p(t) & :=& \gamma_p \langle t \rangle,\qquad   t
\in[0,1]^d,
\\
\rho_p(y) & :=& y^{p} - \gamma_p, \qquad  y
\in\mathbb{R},
\end{eqnarray*}
where $\gamma_p$ is the $p$th moment of the standard Gaussian law,
that is,
\begin{eqnarray*}
&& \gamma_p :=\int_{\mathbb{R}} y^p \gamma(
\mathrm{d}y) =
\cases{ 0, & \quad\mbox{$p$ is odd,}
\vspace*{3pt}\cr
\displaystyle\prod
_{j=1}^{p/2} (2j-1), &\quad \mbox{$p$ is even.}}
\end{eqnarray*}
Since the function $\rho_p$ is a polynomial, it belongs to
$L^2(\mathbb{R}
,\gamma)$ and is a linear combination of finitely many Hermite polynomials.
Moreover, it is easy to check that the Hermite rank of $\rho_p$ is
given by
\begin{eqnarray*}
\underline{k} = \underline{k}_p = %
\cases{ 1, &
\quad\mbox{$p$ is odd,}
\vspace*{3pt}\cr
2, & \quad\mbox{$p$ is even.}}
\end{eqnarray*}
Thus, the Hermite coefficients of $\rho_p$ satisfy Assumption~\ref
{Hermite-assumption}. In what follows, we denote by $\Lambda_{H,\rho
_p}$ the constant given by \eqref{lambdadef}, substituting $f$ with
$\rho_p$ therein.

As a straightforward application of Theorems \ref{gvar-clt} and \ref
{gvar-nclt}, we can prove a functional law of large numbers (FLLN) for
$V^{(n)}_{p}$ as $n \rightarrow\infty$, namely,
\begin{eqnarray*}
&& V^{(n)}_p \mathop{\longrightarrow}_{n\rightarrow\infty}^{\mathbf{P}}
v_p  \qquad \mbox{in }D\bigl([0,1]^d\bigr).
\end{eqnarray*}
It would then be natural to expect that the rescaled fluctuation process
\begin{equation}
\label{pv-fluctuation}
\frac{\langle m(n) \rangle}{\langle c(n) \rangle^{{1}/{2}}} \bigl(V^{(n)}_p(t) -
v_p(t) \bigr), \qquad  t \in[0,1]^d,
\end{equation}
has a non-trivial limit as $n \rightarrow\infty$.
In fact, we can write for any $t \in[0,1]^d$ and $n \in\mathbb{N}$,
\begin{equation}
\label{pv-decomp} \frac{\langle m(n) \rangle}{\langle c(n) \rangle^{{1}/{2}}} \bigl(V^{(n)}_p(t) -
v_p(t) \bigr) = \overline{U}^{(n)}_{\rho_p}(t) -
\beta ^{(n)}_p(t),
\end{equation}
where
\begin{eqnarray*}
&& \beta^{(n)}_p(t) :=\frac{\langle m(n)\rangle}{\langle
c(n)\rangle^{{1}/{2}}}
\biggl(v_p(t) - v_p \biggl( \frac{\lfloor m(n) t
\rfloor
}{m(n)} \biggr)
\biggr)\geq0.
\end{eqnarray*}
If the remainder $\beta^{(n)}_p$ were asymptotically negligible in
$D([0,1]^d)$, the limit of the fluctuation process \eqref
{pv-fluctuation} when $n \rightarrow\infty$ would be easy to deduce
from Theorems \ref{gvar-clt} and \ref{gvar-nclt}.
If $p$ is odd, then indeed $\beta^{(n)}_p=0=v_p$ for any $n \in
\mathbb{N}$.
However, when $p$ is even, the situation is more delicate.
In the special case $d=1$, it is not difficult to see that $\beta
^{(n)}_p(t) < c(n)^{-{1}/{2}} \rightarrow0$ when $n \rightarrow
\infty$ for any $t \in[0,1]$. But when $d\geq2$, the fluctuations of
$\beta^{(n)}_p$ may be non-negligible or even explosive when $n
\rightarrow\infty$, as the following example shows.

\begin{exm}
Consider the case where $p$ is even, $d\geq2$, $m(n) :=(n,\ldots
,n)$ for any $n \in\mathbb{N}$, and $H\in (0,\frac{3}{4}
)^d$. Then we
have by the mean value theorem,
\begin{eqnarray*}
&& \beta^{(n)}_p(t) = n^{d/2-1} \sum
_{\nu=1}^d \biggl( \prod
_{\kappa
\neq
\nu} \xi^{(n)}_\kappa(t) \biggr) \{n
t_\nu\}, \qquad  t \in[0,1]^d,   n \in\mathbb{N},
\end{eqnarray*}
where $\xi^{(n)}(t)$ is some convex combination of $n^{-1}\lfloor
nt\rfloor$ and $t$.
We will now show that $\beta^{(n)}_p$ cannot converge to a continuous
function in $D([0,1]^d)$ as $n \rightarrow\infty$ (similar, but
slightly longer, argument shows that a discontinuous limit in
$D([0,1]^d)$ is also impossible).

To this end, suppose that $\beta^{(n)}_p \rightarrow\beta$ in
$D([0,1]^d)$, where $\beta\in C([0,1]^d)$. Then it follows that $\beta
^{(n)}_p \rightarrow\beta$ uniformly. By the continuity of $\beta$,
there exists an open set $E \subset [\frac{2}{3},1 ]^d$ such that
\begin{equation}
\label{modulus}
\sup_{s,t \in E} \bigl|\beta(s) - \beta(t)\bigr| \leq
\frac{1}{2^{d}}.
\end{equation}
Note that there exists $n_0 \in\mathbb{N}$ such that $E \cap\mathscr{E}_n
\neq\varnothing$ for any $n \geq n_0$, where $\mathscr{E}_n = \{
i/m(n) \dvt 0 \leq i \leq m(n) \}$. Moreover, we can find $n_1 \geq n_0$
such that
\begin{eqnarray*}
&& \inf_{t \in E}\prod_{\kappa\neq\nu}
\xi^{(n)}_\kappa(t) \geq \frac
{1}{2^{d-1}} \qquad \mbox{for
any $n \geq n_1$.}
\end{eqnarray*}
Thus, we find that for any $n \geq n_1$,
\begin{equation}
\label{E-upper} \sup_{t \in E} \beta^{(n)}_p(t)
\geq\frac{n^{d/2-1}}{2^{d-1}},
\end{equation}
while
\begin{equation}
\label{E-lower} \inf_{t \in E} \beta^{(n)}_p(t)
= 0.
\end{equation}
But when $\beta^{(n)}_p \rightarrow\beta$ uniformly, the estimate
\eqref{modulus} is not compatible with \eqref{E-upper} and \eqref
{E-lower}, which is a contradiction. (This also shows that $\beta
^{(n)}_p$ cannot converge to $\beta$ along a subsequence.)
\end{exm}

\subsection{Multilinear interpolations}

We have just seen that the rescaled fluctuations \eqref{pv-fluctuation}
of the power variations $V^{(n)}_p$, $n\in\mathbb{N}$, around their
FLLN limit
$v_p$ do not necessarily satisfy a functional limit theorem in
$D([0,1]^d)$ when $d\geq2$ and $p$ is even. Note that it is implicit in
the definition of $V^{(n)}_p$ that the corresponding partial sums are
interpolated in a \emph{piecewise constant} manner. Such an
interpolation can have very poor precision in higher dimensions. In
fact, interpolating $V^{(n)}_p$ using a more appropriate, multilinear
method enables functional convergence in the general case.

\begin{defn}
For any $n\in\mathbb{N}$, we define a (piecewise) \emph{multilinear
interpolation operator}
$L_{n} \dvtx \mathbb{R}^{[0,1]^d}\rightarrow C([0,1]^d)$ acting on a
function $g \dvtx
[0,1]^d \rightarrow\mathbb{R}$, sampled on the lattice $\mathscr{E}_{n}$,
by
\begin{equation}
\label{interpolation-def}
(L_{n}g) (t) :=\sum_{i \in\{0,1\}^d}
g \biggl(\frac{\lfloor
m(n)t\rfloor+i}{m(n)} \biggr)\alpha^{(n)}_i(t),
\qquad  t \in[0,1]^d,
\end{equation}
where the weights
\begin{eqnarray*}
&&\alpha^{(n)}_i(t) :=\bigl\langle\bigl\{m(n)t\bigr
\}^i \bigl(1-\bigl\{ m(n)t\bigr\} \bigr)^{1-i}\bigr\rangle,
 \qquad i \in\{0,1\}^d,
\end{eqnarray*}
belong to $[0,1]$ and satisfy
\begin{equation}
\label{convex}
\sum_{i \in\{0,1\}^d}\alpha^{(n)}_i(t)
= 1.
\end{equation}
\end{defn}

\begin{rem}
(1)
In the cases $d=1$ and $d=2$, the definition \eqref
{interpolation-def} reduces to the well-known (piecewise) linear and
bilinear interpolation formulae, respectively.

(2) The definition \eqref{interpolation-def} involves slight abuse of
notation. Namely,
\begin{equation}
\label{out-of-bounds}
\frac{\lfloor m(n)t \rfloor+ i}{m(n)} \notin[0,1]^d
\end{equation}
when $t_\nu=1$ and $i_\nu= 1$ for some $\nu\in\{1,\ldots,d\}$. But
then $\alpha^{(n)}_i(t) = 0$, whence \eqref{out-of-bounds} is of no concern.
\end{rem}

The fluctuation process, analogous to \eqref{pv-fluctuation}, obtained
by substituting the power variation $V^{(n)}_p$ with its multilinear
interpolation $\widetilde{V}^{(n)}_p :=L_{n} V^{(n)}_p$ satisfies
the following functional limit theorem. In particular, it applies with
any $d \in\mathbb{N}$ and $p\in\mathbb{N}$.

\begin{thm}[(Interpolated power variations)]
\label{inter-flt}
(1) If $H \in(0,1)^d \setminus (1-\frac{1}{2
\underline{k}_p},1 )^d$, then
\begin{eqnarray*}
&& \biggl(Z, \frac{\langle m(n) \rangle}{\langle c(n) \rangle^{1/2}} \bigl(\widetilde{V}^{(n)}_p-v_p
\bigr) \biggr) \mathop{\longrightarrow}_{n\rightarrow \infty }^{\mathscr{L}} (Z,
\Lambda_{H,\rho_p}\widetilde{Z} ) \qquad \mbox{in $C\bigl([0,1]^d
\bigr)^2$,}
\end{eqnarray*}
where $\widetilde{Z}$ is the fBs\vspace*{1pt} of Theorem~\ref{gvar-clt}.

(2) If $H \in (1-\frac{1}{2 \underline
{k}_p},1 )^d$, then
\begin{eqnarray*}
&& \frac{\langle m(n) \rangle}{\langle c(n) \rangle^{1/2}} \bigl(\widetilde {V}^{(n)}_p-v_p
\bigr) \mathop{\longrightarrow}_{n\rightarrow\infty}^{\mathscr {\mathbf{P} }}
\Lambda_{H,\rho_p} \widehat{Z} \qquad \mbox{in $C\bigl([0,1]^d
\bigr)$,}
\end{eqnarray*}
where $\widehat{Z}$ is the Hermite sheet of Theorem~\ref{gvar-nclt}.
\end{thm}

\begin{rem}
As mentioned above, the remainder term $\beta^{(n)}_p$ in the
decomposition \eqref{pv-decomp} is asymptotically negligible in
$D([0,1]^d)$ if $d=1$ or $p$ is odd. In these special cases,
multilinear interpolations can be dispensed with,
to wit the convergences of Theorem~\ref{inter-flt} hold also with the
original power variation $V^{(n)}_p$ in place of $\widetilde
{V}^{(n)}_p$, in the spaces $D([0,1]^d)^2$ and $D([0,1]^d)$, respectively.
\end{rem}

The proof of Theorem~\ref{inter-flt} is based on the following two
simple lemmas concerning the multilinear interpolation operators.
First, we show that the function $v_p$ is a fixed point of the operator
$L_n$ for any $n \in\mathbb{N}$.

\begin{lem}[(Fixed point)]\label{fixed-point}
We have $L_n v_p = v_p$ for any $n \in\mathbb{N}$.
\end{lem}

\begin{pf}
Let $t\in[0,1]^d$ and $n \in\mathbb{N}$. By rearranging, we obtain that
\begin{eqnarray*}
(L_n v_p) (t) & =& \sum
_{i \in\{0,1\}^d} \gamma_p \biggl\langle
\frac
{\lfloor m(n)t\rfloor+ i}{m(n)} \bigl\{m(n)t\bigr\}^i \bigl(1-\bigl\{ m(n)t\bigr
\}\bigr)^{1-i} \biggr\rangle
\\
& =& \gamma_p \prod_{\nu= 1}^d
\sum_{j\in\{0,1\}} \frac{\lfloor
m_\nu
(n) t_\nu\rfloor+j}{m_\nu(n)}\bigl\{
m_\nu(n) t_\nu\bigr\}^j \bigl(1-\bigl\{
m_\nu(n) t_\nu\bigr\}\bigr)^{1-j}.
\end{eqnarray*}
It remains to observe that for any $\nu\in\{1,\ldots,d\}$,
\begin{eqnarray*}
&&\sum_{j\in\{0,1\}} \frac{\lfloor m_\nu(n) t_\nu\rfloor+j}{m_\nu
(n)}\bigl\{
m_\nu(n) t_\nu\bigr\}^j \bigl(1-\bigl\{
m_\nu(n) t_\nu\bigr\}\bigr)^{1-j} =
\frac{\lfloor
m_\nu
(n) t_\nu\rfloor+ \{ m_\nu(n) t_\nu\}}{m_\nu(n)} = t_\nu,
\end{eqnarray*}
and the assertion follows.
\end{pf}

Second, we show that convergence in probability in the space
$D([0,1]^d)$ can be converted to convergence in probability in
$C([0,1]^d)$ via interpolations.

\begin{lem}[(Convergence and interpolation)]\label{inter-conv}
Let $X_1,X_2,\ldots$ be random elements in $D([0, 1]^d)$ and $X$ a
random element in $C([0,1]^d)$, all defined on a common probability
space. If $X_n \displaystyle\mathop{\longrightarrow}^{\mathbf{P}} X$ in $D([0,1]^d)$ as $n
\rightarrow\infty$, then
\begin{eqnarray*}
&& L_n X_n \mathop{\longrightarrow}_{n \rightarrow\infty}^{\mathbf{P}} X
 \qquad \mbox{in $C\bigl([0,1]^d\bigr)$.}
\end{eqnarray*}
\end{lem}

\begin{pf}
By \eqref{convex}, we can write for any $t \in[0,1]^d$ and $n \in
\mathbb{N}$,
\begin{eqnarray*}
(L_n X_n) (t)-X(t) & =& \sum
_{i \in\{0,1\}^d} \biggl(X_n \biggl(
\frac
{\lfloor m(n)t\rfloor+i}{m(n)} \biggr) - X \biggl(\frac{\lfloor
m(n)t\rfloor
+i}{m(n)} \biggr) \biggr)
\alpha^{(n)}_i(t)
\\
&&{}  + \sum_{i \in\{0,1\}^d} \biggl(X \biggl(
\frac{\lfloor
m(n)t\rfloor+i}{m(n)} \biggr) - X(t) \biggr)\alpha^{(n)}_i(t).
\end{eqnarray*}
Thus, invoking \eqref{convex} again, we obtain the bound
\begin{eqnarray*}
&& \sup_{t \in[0,1]^d}\bigl|(L_n X_n) (t)-X(t)\bigr| \leq
\sup_{t \in
[0,1]^d}\bigl|X_n(t)-X(t)\bigr| + w_X \bigl(
\underline{m}(n)^{-1} \bigr),
\end{eqnarray*}
where
\begin{eqnarray*}
&& w_X(u) :=\sup \bigl\{\bigl|X(s)-X(t)\bigr| \dvt s, t \in[0,1]^d,
\|s-t\| _\infty\leq u \bigr\},\qquad  u > 0,
\end{eqnarray*}
is the modulus of continuity of $X$, which satisfies $\lim_{u
\rightarrow0}w_X(u) = 0$ a.s. since the realizations of $X$ are
uniformly continuous. Thus, $\lim_{n \rightarrow\infty}w_X
(\underline{m}(n)^{-1} ) = 0$ a.s. Finally, since convergence to a
continuous function in $D([0,1]^d)$ is equivalent to uniform
convergence, it follows that $\sup_{t \in[0,1]^d}|X_n(t)-X(t)|\stackrel{\mathbf{P}}{\longrightarrow} 0$ as $n \rightarrow\infty$.
\end{pf}

\begin{pf*}{Proof of Theorem~\ref{inter-flt}}
We have for any $n \in\mathbb{N}$, by Lemma~\ref{fixed-point}, decomposition
\eqref{pv-decomp}, and the linearity of the operator $L_n$,
\begin{eqnarray*}
&& \frac{\langle m(n) \rangle}{\langle c(n) \rangle^{1/2}} \bigl(\widetilde {V}^{(n)}_p-v_p
\bigr) = L_n \biggl( \frac{\langle m(n) \rangle
}{\langle
c(n) \rangle^{1/2}} \bigl(V^{(n)}_p-v_p
\bigr) \biggr) = L_n \overline {U}^{(n)}_{\rho_p} +
L_n \beta_{p}^{(n)}.
\end{eqnarray*}
Note that the function
\begin{eqnarray*}
&& t \mapsto v_p \biggl(\frac{\lfloor m(n) t \rfloor}{m(n)} \biggr)
\end{eqnarray*}
coincides with $v_p$ on $\mathscr{E}_n$. Since $L_n g$ depends on the
function $g$ only through the values of $g$ on $\mathscr{E}_n$, we
find that
\begin{eqnarray*}
&& L_n v_p = L_n v_p \biggl(
\frac{\lfloor m(n) \cdot\rfloor}{m(n)} \biggr),
\end{eqnarray*}
whence
\begin{eqnarray*}
&& L_n \beta_{p}^{(n)} = \frac{\langle m(n) \rangle}{\langle c(n)
\rangle
^{1/2}}
\biggl( L_n v_p - L_n v_p
\biggl(\frac{\lfloor m(n) \cdot
\rfloor
}{m(n)} \biggr) \biggr) = 0.
\end{eqnarray*}
The assertion in the case (2) follows now from Theorem~\ref{gvar-nclt} and Lemma~\ref{inter-conv}. In the case
(1), one can apply Theorem~\ref{gvar-clt}, Lemma~\ref
{inter-conv} and Skorohod's representation theorem \cite{Kal2002}, Theorem~4.30.
\end{pf*}

\section*{Acknowledgements}

M.S. Pakkanen wishes to thank CEREMADE, Universit\'e Paris-Dauphine
for warm hospitality and acknowledge support
from CREATES (DNRF78), funded by the Danish National Research Foundation,
from the Aarhus University Research Foundation regarding the project
``Stochastic and Econometric Analysis of Commodity Markets,'' and
from the Academy of Finland (project 258042).

% imsref loaded by daiva.urboniene, 2015-05-05 09:26:53

\printhistory
\end{document}